\documentclass[11pt,amssymb]{amsart}
\usepackage{amsmath}
\usepackage{amssymb}
\DeclareFontFamily{OT1}{rsfs}{}
\DeclareFontShape{OT1}{rsfs}{n}{it}{<-> rsfs10}{}
\DeclareMathAlphabet{\mathscr}{OT1}{rsfs}{n}{it}
\def\ni{\noindent}
\def\bA{\mathbf{A}}
\def\bB{\mathbf{B}}
\def\bG{\mathbf{G}}
\def\bnG{\mathbf{G}^{\natural}}
\def\bnP{\mathbf{P}^{\natural}}
\def\bnQ{\mathbf{Q}^{\natural}}
\def\bC{\mathbf{C}}
\def\bF{\mathbf{F}}
\def\bH{\mathbf{H}}
\def\bN{\mathbf{N}}
\def\bM{\mathbf{M}}
\def\bR{\mathbf{R}}
\def\bS{\mathbf{S}}
\def\bT{\mathbf{T}}
\def\bP{\mathbf{P}}
\def\bQ{\mathbf{Q}}
\def\bU{\mathbf{U}}
\def\bV{\mathbf{V}}

\def\sX{\mathscr{X}}
\def\sN{\mathscr{N}}
\def\sB{\mathscr{B}}
\def\sD{\mathscr{D}}
\def\sG{\mathscr{G}}
\def\sH{\mathscr{H}}
\def\sN{\mathscr{N}}
\def\sP{\mathscr{P}}
\def\sR{\mathscr{R}}
\def\sS{\mathscr{S}}
\def\sU{\mathscr{U}}

\def\sZ{\mathscr{Z}}
\def\ok{\overline{k}}

\def\oH{\overline{H}}
\def\oU{\overline{U}}

\def\snG{\mathscr{G}^{\natural}}

\def\osM{\overline{\mathscr{M}}}

\begin{document}

\vskip7mm

\centerline{\bf Weakly-split spherical Tits systems in quasi-reductive groups}
\vskip5mm

\centerline{By {\sc Gopal Prasad}}

\vskip6mm

 \centerline{\it Dedicated to C.\,S.\,Seshadri on his 80th birthday}\vskip5mm

 {\it Abstract.} We will prove that any weakly-split spherical Tits system $(B,N)$ in $G= \bG(k)$ ($\bG$ a quasi-reductive $k$-group, such as a connected reductive  $k$-group) satisfying 
 some natural conditions is ``standard''. In particular, if $\bG$ is anisotropic over $k$, then such a  Tits system is trivial, i.e., $B= G$.
\vskip4mm
   
\centerline{\bf 1. Introduction}

\vskip3mm

\ni{\bf 1.1.} In this paper, $k$ will always be an infinite field, $\ok$  a fixed algebraic closure of $k$ and $k_s$ the separable closure of $k$ in $\ok$.  
\vskip1mm

It is known that all maximal $k$-split tori in a smooth connected affine $k$-group $\bG$ are $\bG(k)$-conjugate ([6], Theorem C.2.3). The {\it $k$-rank} of such a $\bG$ is by definition the dimension of a maximal $k$-split torus of $\bG$. We will say that a smooth connected affine $k$-group is $k$-{\it isotropic} if its $k$-rank is positive, and $k$-{\it anisotropic} if its $k$-rank is 0.  
\vskip1mm

A smooth connected unipotent $k$-group $\bU$ is said to be {\it $k$-wound} if every map of $k$-schemes $\bA^1_k\rightarrow \bU$ is a constant map to a point in $\bU(k)$.   It is obvious from this definition that  such a group does not contain any nontrivial $k$-split smooth connected $k$-subgroup.    It is known that  a $k$-torus can only act trivially on a $k$-wound $\bU$  (see [6], Proposition B.4.4). Therefore, if $\bU$ is a $k$-wound smooth connected unipotent normal $k$-subgroup of a $k$-group $\bH$, then every $k$-torus of $\bH$ commutes with $\bU$.
According to Theorem B.3.4 of [6], any smooth connected unipotent $k$-group $\bU$ contains a unique $k$-split smooth connected normal $k$-subgroup, to be  denote by ${\bU}_{\rm{split}}$, such that $\bU/\bU_{\rm{split}}$ 
is $k$-wound. Moreover,  this subgroup contains every  $k$-split smooth connected $k$-subgroup of $\bU$, so it can alternatively be described as the unique maximal $k$-split  smooth connected $k$-subgroup of $\bU$.  Therefore, a smooth connected unipotent $k$-group  is $k$-wound if  it does not contain any nontrivial $k$-split smooth connected $k$-subgroup.
\vskip1mm

The maximal smooth connected unipotent normal $k$-subgroup of $\bG$ is called the $k$-{\it unipotent radical}  of  $\bG$; it will be denoted 
by $\sR_{u,k}(\bG)$. The maximal $k$-split  smooth connected unipotent normal $k$-subgroup of  $\bG$  is called the $k$-{\it split unipotent radical} of $\bG$; it will be denoted by $\sR_{us,k}(\bG)$, and is clearly contained in $\sR_{u,k}(\bG)$. By Corollary B.3.5 of [6], the $k$-split unipotent radical $\sR_{us,k}(\bG)$ is the maximal smooth connected $k$-split subgroup of $\sR_{u,k}(\bG)$, so the quotient $\sR_{u,k}(\bG)/\sR_{us,k}(\bG)$ is $k$-wound.
\vskip1mm

A smooth connected affine $k$ group $\bG$ is said to be {\it pseudo-reductive} if its $k$-unipotent radical $\sR_{u,k}(\bG)$ is trivial, and is said to be {\it quasi-reductive} if its $k$-split unipotent radical $\sR_{us,k}(\bG)$ is trivial. If $\bG$ is quasi-reductive then its $k$-unipotent radical $\sR_{u,k}(\bG)$ equals $\sR_{u,k}(\bG)/\sR_{us,k}(\bG)$, and hence is $k$-wound. 
Note that the notion of a quasi-reductive group is different from the notion of a quasi-reductive group scheme over a discrete valuation ring introduced in a joint paper of the author with Jiu-Kang Yu (in J.\,Alg.\,Geom.\,{\bf 15}(2006)).   
\vskip1mm

Given a smooth connected affine $k$-group $\bG$,  $\bG/\sR_{u,k}(\bG)$ and $\bG/\sR_{us,k}(\bG)$ are respectively the maximal pseudo-reductive and the maximal  quasi-reductive quotients of $\bG$. If $k$ is perfect, every pseudo-reductive and every quasi-reductive $k$-group is reductive. 
\vskip1mm

In the rest of the paper, $\bG$ will be a quasi-reductive $k$-group, $\sD(\bG)$ will denote its derived subgroup $(\bG,\bG)$.  Being $k$-wound, $\sR_{u,k}(\bG)$ commutes with every $k$-torus of $\bG$ (Proposition B.4.4 of [6]). In particular,  if $\bS$ is a maximal $k$-split torus of $\bG$ then $\sR_{u,k}(\bG)\subset Z_{\bG}(\bS)$, so the $k$-root system (see C.2.12 in [6]) of $\bG$ with respect to $\bS$ is the set of nonzero weights of $\bS$ on ${\rm{Lie}}(\bG)$.    There is a natural bijective correspondence between the 
set of $k$-isotropic  minimal perfect smooth connected proper normal $k$-subgroups of $\bG$ and the set of irreducible components of the $k$-root system of $\bG$ with respect to $\bS$; see Proposition C.2.32  in [6].
\vskip1mm

\vskip1mm

\ni{\bf 1.2. The standard Tits systems in $G:=\bG(k)$ and in certain subgroups of $G$.}  Let $G = \bG(k)$ and $G^+$ be the normal subgroup of $G$ generated by the group of $k$-rational points of the $k$-split unipotent radicals  of pseudo-parabolic $k$-subgroups of $\bG$ (for definition of pseudo-parabolic $k$-subgroups, see [6], 2.2; in a connected reductive group the pseudo-parabolic $k$-subgroups 
are the parabolic $k$-subgroups by Proposition 2.2.9 of [6]).  If $\bG$ is $k$-isotropic and it does not contain a nontrivial perfect smooth connected proper normal $k$-subgroup, then $G^+$, and in fact every non-central normal subgroup of $G$, is Zariski-dense in $\bG$ (Theorem C.2.33  of [6]).   \vskip1mm

Let $\bS$ be a maximal $k$-split torus of $\bG$. Let  $Z_{\bG}(\bS)$ and $N_{\bG}(\bS)$ respectively be the centralizer and the normalizer of $\bS$ in $\bG$. The finite group $_kW :=N_{\bG}(\bS)(k)/Z_{\bG}(\bS)(k)$ is called the {\it $k$-Weyl group} of $\bG$, it is the Weyl group of the $k$-root system $_k\Phi$ of $\bG$ with respect to $\bS$. 

For simplicity, we will denote the $k$-root system $_k\Phi$ of $\bG$ and the $k$-Weyl group $_kW$ by $\Phi$ and $W$ respectively in the sequel.  
\vskip1mm

Let $\bP$ be a minimal pseudo-parabolic $k$-subgroup  of $\bG$ containing $\bS$. Then $\bP$ contains the centralizer $Z_{\bG}(\bS)$ of $\bS$ in $\bG$ ([6], Proposition C.2.4). 
We fix a subgroup $\sG$ of $G$ containing $G^+$ such that the identity component of the Zariski-closure of $\sG$ in $\bG$ is same as the identity component of the Zariski-closure of $G$.  Let $\sB = \bP(k)\cap \sG$ and $\sN = N_{\bG}(\bS)(k)\cap \sG$.    It follows from Proposition C.2.23 of [6] that $\sN$ maps onto the Weyl group $W$. It has been observed in Remark C.2.24 of [6] (using Theorem C.2.19 and Remark C.2.21 of [6]) that 
$(\sB,\sN)$ is a Tits system in $\sG$. The Weyl group of this Tits system is $W$ and the rank of this Tits system is equal to the rank of the $k$-root system $\Phi$ of $\bG$ (which is equal to the $k$-rank of $\sD(\bG)$;  see Theorem C.2.14 
of [6]).   We call this Tits system a  {\it standard Tits system in $\sG$}. (Note that this notion of a standard Tits system in the abstract group $\sG$ makes use of  the given quasi-reductive 
algebraic $k$-group $\bG$ whose group of $k$-rational points contains the group $\sG$. But a consequence of Theorem B of this paper  is that  the notion is to a large extent independent of $\bG$.) Conjugacy of maximal $k$-split tori and of minimal pseudo-parabolic $k$-subgroups  in $\bG$ under $\sG$ (see Remark C.2.24 of [6]) implies that any two standard Tits systems in $\sG$ are conjugate to each other.
\vskip1mm

The Bruhat decomposition for the Tits system $(\sB,\sN)$ of $\sG$ gives us that $\sG$ is a disjoint union of $\sB w\sB$, $w\in W$.
\vskip1mm

It turns out, rather surprisingly, that any Tits system in $\sG$ satisfying some natural conditions is a standard Tits system. For a precise statement see Theorem B below.  

\vskip2mm
\ni{\bf 1.3.} From Remark C.2.11 and Lemma 1.2.5 of [6] applied to a maximal $k$-torus of $\bG$ containing the given maximal $k$-split torus $\bS$ of $\bG$, we see that $\bS$ is an almost direct product of the maximal $k$-split central torus of $\bG$ and  the maximal $k$-split torus of $\sD(\bG)$ contained in $\bS$. Thus, $\sD(\bG)$ is $k$-anisotropic if and only if $\bS$ is central. By conjugacy of maximal $k$-split tori in $\bG$, the centrality of $\bS$ is equivalent to the assertion that every $k$-split torus of $\bG$ is central.  Since, as mentioned above, $\sR_{u,k}(\bG)$ commutes with every $k$-torus of $\bG$, we easily see that the condition that every $k$-split torus of $\bG$ is central is equivalent to the condition that every $k$-split torus of the maximal pseudo-reductive quotient $\bG/\sR_{u,k}(\bG)$ is central. But according to Lemma 2.2.3 of [6], this latter condition is equivalent to the condition that $\bG$ does not contain proper pseudo-parabolic $k$-subgroups.   Hence if every $k$-split torus of $\bG$ is central (or, equivalently, if $\sD(\bG)$ is $k$-anisotropic), then $G^+$ is trivial.

\vskip1mm

\ni{\bf 1.4.} Let $(B,N)$ be an arbitrary Tits system in $\sG$. Let $H= B\cap N$. Then the {\it Weyl group} of this Tits system is $N/H$; we will denote it by $W^T$ in the sequel.  For basic results on Tits systems, including those recalled in 3.1 below, see [3], Ch.\,IV.  

Let $S$ be the distinguished  set of involutive generators of  $W^T$. The {\it rank} of the Tits system $(B,N)$ is the cardinality of $S$. For $s\in S$, let $\sG_s = B\cup BsB$. Then $\sG_s$ is a subgroup of $\sG$ for every $s\in S$.

We will say that the Tits system $(B,N)$ is {\it weakly-split} if there exists a nilpotent normal subgroup $U$ of $B$ such that  
$B = HU$.
We will say that the Tits system  is {\it split}\,  if it is {\it saturated} (that is, $H =\bigcap_{n\in N}nBn^{-1}$) and there exists a nilpotent normal subgroup $U$ of $B$ such that $B = H\ltimes U$. As $\sR_{u,k}(\bG)\subset Z_{\bG}(\bS)$, the standard Tits system in $\sG$ with $\sB =\bP(k)\cap \sG$ and $\sN =N_{\bG}(\bS)(k)\cap \sG$ is split, with $\sB = \sH\ltimes \sU$, where $\sH = \sB\cap \sN = Z_{\bG}(\bS)(k)\cap \sG$ and $\sU =\sR_{us,k}(\bP)(k)$ (see Remarks C.2.21 and C.2.24  in [6]).

\vskip2mm

In this paper, we will only work with {\it spherical} Tits systems  $(B,N)$ in $\sG$ (i.e., Tits systems with finite Weyl group). The Tits systems being considered here will often be either weakly-split  or split, and $U$ will be a nilpotent normal subgroup of $B$ as in 1.4.
\vskip1mm

The purpose of this paper is to prove the following two theorems.
\vskip2mm

\ni {\bf Theorem A.} {\em Assume that  $({\mathrm{i}})$\,either $\bG$ is perfect and quasi-reductive, or it is  pseudo-reductive, $({\mathrm{ii}})$\,every $k$-split torus of $\bG$ is central (or, equivalently, $\sD(\bG)$ is $k$-anisotropic, cf.\,1.3),  and $({\mathrm{iii}})$\,the Tits system $(B,N)$ is weakly-split. Then  $B$ is of finite index in $\sG$.  Furthermore,  if one of the following two conditions hold, then $B=\sG$: 
\vskip1mm
\begin{enumerate}
\item $\bG$ is reductive.
\item $\bG$ is perfect and quasi-reductive. 
\end{enumerate}}

\ni{\bf Theorem B.} {\em  Assume that $\bG$ is an arbitrary quasi-reductive group, the Tits system $(B,N)$ is weakly-split,  and for every $s\in S$, the index of $B$ in $\sG_s$ is infinite (or, equivalently, the index of $B\cap sBs^{-1}$ in $B$ is infinite).  Then there exists a pseudo-parabolic $k$-subgroup $\bP$ of $\bG$ such that $B = \bP(k)\cap \sG$.
\vskip1mm

Assume further that the Tits system $(B,N)$ is saturated and $B$ does not contain a non-central normal subgroup of the group of $k$-rational points of a $k$-isotropic minimal perfect smooth connected normal $k$-subgroup  of $\mathbf{G}$.  Then $(B,N)$ is a standard Tits system,  namely it is one of the Tits systems given in 1.2 above. In particular, the Tits system is split.}

\vskip2mm

Note that for the standard Tits system in $\sG$, the index of $B$ in $\sG_s$ is infinite for every $s\in S$ (since the field $k$ has been assumed to be infinite). The standard Tits system is saturated and split (1.4). 
\vskip2mm

\ni{\it Remark 1.}  Regarding the condition in Theorem B that ``$B$ does not contain a non-central normal subgroup of the group of $k$-rational points of a $k$-isotropic minimal perfect smooth connected normal $k$-subgroup  of 
$\mathbf{G}$'' we make the following observation:  Theorem C.2.32\,(2) of [6]  implies that the intersection of $G^+$ with the group $\bN(k)$  of $k$-rational points of a $k$-isotropic minimal perfect smooth 
connected normal $k$-subgroup $\bN$ of $\bG$ is a normal subgroup of $\bN(k)$ that is not central.  We also mention that $\bN$ does not contain a nontrivial perfect smooth connected proper $k$-subgroup that is normal in $\bN$ (Proposition C.2.32\,(4) of [6]), and so by Theorem C.2.33\,(3) of [6] every non-central normal subgroup of $\bN(k)$ 
is Zariski-dense in $\bN$. 
\vskip1mm

 Let $\bG'$ be a $k$-isotropic $k$-simple connected semi-simple $k$-group and $G' = \bG'(k)$. Let $(B,N)$ be a weakly-split Tits system in $G=\bG(k)$. Then $(B\times G', N\times G')$ is a weakly-split Tits system in $G\times G'$. This shows  that the condition imposed on $B$ in the second assertion of Theorem B is necessary.

\vskip2mm
Theorems A and B hold in particular for $\bG$ connected reductive, and appear to be new and interesting even for these groups.  We recall that there are spherical Moufang buildings which arise from nonreductive pseudo-reductive groups (see, for example, [13], 10.3.2, or [14], 41.20; the groups appearing in 10.3.2 of [13] are the group of $k$-rational points of ``exotic'' pseudo-reductive groups described in [6], Ch.\,7). Therefore, we have chosen to work in the more general set-up of quasi-reductive groups in this paper.  Proofs of Theorems A and B given below for $\bG$ reductive do not require familiarity with the theory of pseudo-reductive and quasi-reductive groups. It would be enough to know  the Borel-Tits theory of reductive groups over non-algebraically closed fields (as described  in [2]) and to know that the pseudo-parabolic $k$-subgroups of a connected reductive $k$-group are just the parabolic $k$-subgroups (Proposition 2.2.9 of [6]). 
\vskip1mm

We will prove Theorem A in \S4 and prove Theorem B, which is the main result of this paper, in \S5 for Tits systems of rank 1 and in \S6 for Tits systems of arbitrary rank.
\vskip1.5mm

Special cases of Theorem A were proven earlier in [5] and  [1].  Caprace and Marquis proved in [5] that if $\bG$ is semi-simple and $k$-anisotropic, then $B$ is of finite index in $G$ in case $k$ is either perfect or it is a nondiscrete locally compact field. That paper inspired us to write this paper. In [1], Abramenko and Zaremsky have proved Theorem A, by an entirely different argument, for some  classes of semi-simple groups.

\vskip2mm

We owe the following remark to Richard Weiss.
\vskip1mm

\ni{\it Remark 2.} Let $(B,N)$ be a weakly-split spherical Tits system in an abstract group $G$, and let $\Omega$ be the building associated to this Tits system. We assume that (i) $B$ acts faithfully on $\Omega$,  (ii) the Weyl group of the Tits system is irreducible, and (iii) the rank of the Tits system is at least 2. Let $B =HU$, with $H = B\cap N$ and $U$ a nilpotent normal subgroup of $B$. Then the building $\Omega$ is Moufang and $U$ is precisely the subgroup generated by the ``root groups'', of the group ${\rm {Aut}}(\Omega)$ of type-preserving automorphisms of $\Omega$, contained in $B$. This has been proved by De Medts, Haot, Tent and Van Maldeghem, see the corollary in \S2 of [11]  (and also [9] and [10]). 
\vskip1mm

In the early 1970s, Tits proved that every irreducible thick spherical building of rank $\geqslant 3$ is Moufang; see [13], Addenda. Given a semi-simple $k$-group $\bG$, in [13], 11.14, there is an example of a spherical Tits system  $(B,N)$ in a free group $F$ with infinitely many generators such that $B$ does not contain any nontrivial normal subgroup of $F$ and the  building associated to this Tits system  is the building of the standard Tits system in $\bG(k)$. As subgroups of $F$ are free, it is obvious that no Tits system in $F$ can be weakly-split. Thus spherical Moufang buildings can arise from non-weakly-split Tits systems.
\vskip1mm

\ni{\it Remark 3.} Theorem B does not hold in general if $k$ is a finite field and $\bG$ is a connected reductive $k$-group. To  see some examples, let $\bF_p$ denote the finite field with $p$ elements. As mentioned by Katrin Tent, there is a well-known  isomorphism between ${\mathrm{PGL}}_3(\bF_2)$ and ${\mathrm{PGL}}_2(\bF_7)^+$  which provides split Tits systems  in ${\mathrm{PGL}}_3(\bF_2)$ of ranks 1 and 2. There is also an isomorphism between ${\mathrm {PU}}_4(\bF_2)$ and ${\mathrm{PSp}}_4(\bF_3)^+$ which provides two non-isomorphic split Tits systems in ${\mathrm{PU}}_4(\bF_2)$ of rank  2. 

\vskip5mm

\centerline{\bf 2. Preliminaries} 

\vskip3mm

\ni{\bf 2.1.} As in Theorems A and B, $\bG$ is a quasi-reductive $k$-group. Its $k$-unipotent radical $\sR_{u,k}(\bG)$ is $k$-wound, so it commutes with every $k$-torus of $\bG$. 

For a $k$-subgroup $\bH$ of $\bG$, $\bH_{\ok}$ will denote the $\ok$-group obtained from $\bH$ by extension of scalars $k \hookrightarrow \ok$, and $\sR_u(\bH_{\ok})$ will denote the unipotent radical of $\bH_{\ok}$, 
i.e., the maximal smooth connected unipotent normal $\ok$-subgroup of $\bH_{\ok}$. Let $\bG'= \bG_{\ok}/\sR_u(\bG_{\ok})$ be the maximal reductive quotient of $\bG_{\ok}$ and $\pi: \bG_{\ok}\rightarrow \bG'$ be the quotient map. For a $k$-subgroup $\bH$ of $\bG$, we shall denote the image $\pi(\bH_{\ok})\,(\subseteq \bG')$ of $\bH_{\ok}$ by $\bH'$. 

For any smooth affine $k$-group $\bH$, $\sD(\bH)$ will denote the derived subgroup $(\bH,\bH)$ of $\bH$.

\vskip1mm

\ni{\bf 2.2.} For any smooth connected affine $k$-group $\bH$, let $\bH_t$ denote the $k$-subgroup generated by all the $k$-tori of $\bH$. It is known (see Proposition A.2.11 of [6]) that $\bH_t$ is a normal subgroup of $\bH$, $\bH/\bH_t$ is a unipotent group, and  $(\bH_t)_{\ok}$ is the subgroup of $\bH_{\ok}$ generated by all maximal $\ok$-tori of the latter.  Since $\bH/\bH_t$ is a unipotent group, every perfect smooth closed $k$-subgroup of $\bH$ is contained in $\bH_t$, and 
$\bH_t =\bH$ if $\bH$ is perfect.
\vskip1mm

Being generated by $k$-tori, $\bH_t$ is unirational over $k$, and hence $\bH_t(k)$ is Zariski-dense in $\bH_t$.
\vskip1mm

 If $\sR_{u,k}(\bH)$ is $k$-wound, then it commutes with every $k$-torus of $\bH$, and hence it commutes with $\bH_t$.  If, moreover, $\bH$ is perfect,  then as $\bH_t = \bH$,  $\sR_{u,k}(\bH)$ is contained in the center of $\bH$.

\vskip1mm
   
\ni{\bf 2.3. The Zariski-closure of $G = \bG(k)$:} Let $\bnG$ be the identity component of the Zariski-closure of $G$ in $\bG$. Since $\bG_t(k)$ is Zariski-dense in $\bG_t$\,(2.2), $\bG_t\subseteq \bnG$. If either $\bG$ is perfect or reductive, then $\bG_t =\bG$, so in these cases, $\bnG = \bG$; i.e., $\bG(k)$ is Zariski-dense in $\bG$. (Note that  Zariski-density of $G$ in $\bG$ may fail if the latter is a commutative pseudo-reductive group, see Example 11.3.1 in [6].) Since any perfect smooth connected $k$-subgroup of $\bG$ is contained in $\bG_t$, any such subgroup is contained in $\bnG$.

 Since $\bG/\bG_t$ is a unipotent group, under the quotient map $\pi:\bG_{\ok}\rightarrow \bG'$, the subgroup $(\bG_t)_{\ok}$, and so also $\bnG_{\ok}$, maps onto the reductive group $\bG'$. This implies that the geometric unipotent radical $\sR_u(\bnG_{\ok})$ of $\bnG$ is contained in the geometric unipotent radical $\sR_u(\bG_{\ok})$ of $\bG$, which in turn implies that $\sR_{u,k}(\bnG)\subseteq \sR_{u,k}(\bG)$ and $\sR_{us,k}(\bnG)\subseteq \sR_{us,k}(\bG)$. So $\sR_{us,k}(\bnG)$ is trivial, i.e., $\bnG$ is quasi-reductive, and  $\sR_{u,k}(\bnG)$ is $k$-wound.

Given a $k$-torus $\bT$ of $\bG$ and a $k$-split smooth connected unipotent $k$-subgroup $\bU$ of $\bG$, as $\bT(k)$ and $\bU(k)$ are Zariski-dense in $\bT$ and $\bU$ respectively, we see that $\bnG \supset \bT,\, \bU$.  Hence, the root groups and the root systems of $\bG$ and $\bnG$ with respect to any maximal $k$-split torus $\bS$ are the same (see C.2.12 and C.2.20 of [6]), and so the natural homomorphism $$N_{\bnG}(\bS)(k)/Z_{\bnG}(\bS)(k)\longrightarrow N_{\bG}(\bS)(k)/Z_{\bG}(\bS)(k),$$ between the $k$-Weyl groups of $\bnG$ and $\bG$ is an isomorphism.

The group $G^{\natural} := \bnG(k)$ is a normal subgroup of $G$ of finite index, and so it is dense in $\bnG$ in the Zariski-topology.  Recall that we have assumed that the identity component of the Zariski-closure of $\sG$ in $\bG$ is same as the identity component of the Zariski-closure of $G$. So the identity component of the Zariski-closure of $\sG$ in $\bG$ is $\bnG$. The subgroup $\snG:= \bnG(k)\cap \sG$ is a normal subgroup of $\sG$ of finite index and so it is also dense in $\bnG$ in the Zariski-topology. 
\vskip1mm

\ni{\bf 2.4.} For a quasi-reductive group $\bG$,  since $\sR_{u,k}(\bG)\subset Z_{\bG}(\bS)$ for any $k$-torus $\bS$ of $\bG$, a pseudo-parabolic $k$-subgroup $\bP$ of  $\bG$ equals $P_{\bG}(\lambda)$ (in the notation of [6], 2.1 and 2.2, which we will use here and in the sequel) for a $1$-parameter subgroup $\lambda:{\mathrm{GL}}_1\rightarrow \bG$. Since $\lambda({\mathrm{GL}}_1)\subset \bnG$, we see that $\bP^{\natural} := \bP\cap \bnG =P_{\bnG}(\lambda)$ is a pseudo-parabolic $k$-subgroup of $\bnG$. Now let $\bP$ and $\bQ$ be two pseudo-parabolic $k$-subgroups of $\bG$ and $\bnP:= \bP\cap \bnG$, and $\bnQ :=\bQ\cap \bnG$. Assume that $\bnP = \bnQ$. Let $\lambda$ and $\mu$ be $1$-parameter subgroups of $\bG$ such that $\bP = P_{\bG}(\lambda)=Z_{\bG}(\lambda)\ltimes U_{\bG}(\lambda)$ and $\bQ =P_{\bG}(\mu)= Z_{\bG}(\mu)\ltimes U_{\bG}(\mu)$. As the $k$-split smooth connected unipotent groups $U_{\bG}(\lambda)$ and $U_{\bG}(\mu)$ are contained in $\bnG$, $U_{\bnG}(\lambda) = U_{\bG}(\lambda)$ and $U_{\bnG}(\mu) = U_{\bG}(\mu)$. Now since  $\bnP = P_{\bnG}(\lambda)$ and $\bnQ=P_{\bnG}(\mu)$, $\sR_{us,k}(\bnP) = U_{\bG}(\lambda)=\sR_{us,k}(\bP)$ and $\sR_{us,k}(\bnQ)= U_{\bG}(\mu) =\sR_{us,k}(\bQ)$ (Corollary 2.2.5 of [6]). As $\bnP =\bnQ$, we see that $\sR_{us,k}(\bP) = \sR_{us,k}(\bQ)$. Now Corollary 3.5.10 of [6], applied to $K=k_s$, implies that $\bP=\bQ$.  In particular, letting $\bP = \bG$, we see that if $\bnQ = \bnG$, then $\bQ = \bG$. 
\vskip1mm

\ni{\bf 2.5.} Since the multiplication map $$U_{\bnG}(-\lambda)\times Z_{\bnG}(\lambda)\times U_{\bnG}(\lambda) \longrightarrow \bnG$$ is an open immersion (Proposition 3.1.8 of [6]), and the subgroups $U_{\bnG}(-\lambda)(k)$ and $U_{\bnG}(\lambda)(k)$ are contained in $\snG$, the Zariski-density of $\snG$ in 
$\bnG$ implies that $Z_{\bnG}(\lambda)(k)\cap \snG$ and $\bnP(k)\cap \snG$\,($=P_{\bnG}(\lambda)(k)\cap \snG$) are Zariski-dense in $Z_{\bnG}(\lambda)$ and $\bP^{\natural}$\,($=P_{\bnG}(\lambda)= Z_{\bnG}(\lambda)\ltimes U_{\bnG}(\lambda)$) respectively. 
\vskip1mm

Recall that a  $1$-parameter subgroup $\lambda$ of a maximal $k$-split torus $\bS$ of $\bG$ is said to be {\it regular} if for every root $a$ of $\bG$ with respect to $\bS$, $\langle a, \lambda\rangle \ne 0$, or, equivalently, $Z_{\bG}(\lambda) = Z_{\bG}(\bS)$.  In the preceding paragraph, taking $\lambda$ to be a regular $1$-parameter subgroup of a  maximal $k$-split torus $\bS$ of $\bG$, we see that $Z_{\bnG}(\bS)(k)\cap \snG$ is Zariski-dense in $Z_{\bnG}(\bS)$. Hence, the identity component of the Zariski-closure of $Z_{\bG}(\bS)(k)\cap \sG$ is $Z_{\bnG}(\bS)\,(\supseteq\bS)$.  
\vskip1mm

Since $\bnP =\bP\cap \bnG$, we have $\bnP(k)\cap \snG = \bP(k)\cap \snG$. Now the Zariski-density of $\bnP(k)\cap \snG$ in $\bnP$ implies that the identity component of the Zariski-closure of $\bP(k)\cap \snG$, and so also of $\bP(k)$,  is $\bnP$. Therefore, every perfect smooth connected $k$-subgroup, every $k$-torus, and every $k$-split smooth connected unipotent $k$-subgroup of $\bP$ is contained in $\bnP$. If $a$ is a root of $\bG$ with respect to a maximal $k$-split torus $\bS$ contained in $\bP$ and $a$ is a weight of $\bS$ on ${\rm{Lie}}(\bP)$, then as the root group $\bU_a$ is contained in $\bP$, it is also contained in $\bnP$. This implies, in particular, that the weight spaces in ${\rm{Lie}}(\bP)$ and ${\rm{Lie}}(\bnP)$ with respect to $\bS$, for a nonzero weight, are identical.    

\vskip1mm

Let $\bQ$ be a pseudo-parabolic $k$-subgroup of $\bG$ such that $\bQ(k) \supseteq \sG$. Then  $\bnQ :=\bQ\cap \bnG$ equals $\bnG$, and so $\bQ = \bG$ (2.4).

\vskip1.5mm

\ni{\bf 2.6.} We observe here for later use that if  $\bP$ is a pseudo-parabolic $k$-subgroup of  $\bG$ such that $\bP(k)\cap \sG$ is a maximal proper subgroup of $\sG$, then $\bP$ is a maximal proper pseudo-parabolic $k$-subgroup of $\bG$. For, if there is a proper pseudo-parabolic $k$-subgroup $\bQ$ of $\bG$ properly containing $\bP$, then the pseudo-parabolic $k$-subgroup $\bnQ := \bQ\cap \bnG$ of $\bnG$ would properly contain the pseudo-parabolic $k$-subgroup $\bnP := \bP\cap \bnG$ of $\bnG$ (2.4). Now as the identity component of the Zariski-closure of $\bP(k)\cap \sG$ (resp.\,$\bQ(k)\cap \sG$) is $\bnP$ (resp.\,$\bnQ$), we conclude that  $\bP(k)\cap \sG \ne \bQ (k)\cap \sG$.  Hence, $\bQ(k)\supseteq \sG$, which implies that $\bQ = \bG$ (2.5), a contradiction.  

\vskip5mm

\centerline{\bf 3. Six  basic propositions}
\vskip4mm

\ni{\bf 3.1.} As in \S 2, let $\bG$ be a quasi-reductive $k$-group, $G = \bG(k)$. Let $(B,N)$ be a weakly-split Tits system in $\sG$, with $B = HU$, $H = B\cap N$.  Let $W^T = N/H$ be the Weyl group of this Tits system. We assume that the Tits system is spherical, i.e., $W^T$ is finite. Let $S$ be the distinguished  set of involutive generators of $W^T$. For a subset $X$ of $S$, we will denote by $W^T_X$ the subgroup of $W^T$ generated by the elements in $X$, and by $\sG_X$ the subgroup $BW^T_XB$ of  $\sG$.  It is known (as for any Tits system) that $B$ equals its normalizer in $\sG$, so in particular, it  contains the center of $\sG$. (If the Tits system $(B,N)$ is saturated, then the center of $\sG$ is contained in $H$ since $H =\bigcap_{n\in N}nBn^{-1}$ in that case.)  Any subgroup of $\sG$ containing $B$ equals $\sG_X$ for a unique subset $X$ of $S$.  \vskip1mm

\ni{\bf 3.2. Notation.} As in 2.3, we will denote the identity component of the Zariski-closure of $G$ in $\bG$  by $\bnG$. By our assumption, $\bnG$ is also the identity component of the Zariski-closure of $\sG$ in $\bG$. 
Recall that $\bnG$ contains $\bG_t$. Let $\bB$ (resp.,\,$\bU$) be the identity component of the Zariski-closure of   $B$ (resp.,\,$U$) in $\bG$. The subgroups $\bB$ and $\bU$ are contained in 
$\bnG$, and $\bU$ is a connected nilpotent normal $k$-subgroup of $\bB$. Let $B_0=  B\cap \bB(k)$. Then $B_0$ is a normal subgroup of $B$ of finite index. Since $\sG$ is the union of finitely many double cosets 
$BwB$, $w\in W^T$, finitely many $\{B_0,B_0\}$-double cosets cover $\sG$, which implies that some $\{\bB,\bB\}$-double coset  of a $k$-point in $\bnG$ is dense (and hence open) in $\bnG$.  As $\bnG_{\ok}$ maps onto the reductive quotient $\bG'$ of $\bG_{\ok}$ (2.3), we see that some $\{\bB',\bB'\}$-double coset is open in $\bG'$.

\vskip2mm

We shall prove six propositions in this section which will be crucial in this paper.  
\vskip2mm

\ni {\bf Proposition 1.} {\em The following hold:
\vskip1mm

(i) The unique maximal $k$-torus of $\bU$ is contained in the center of $\bG$.
\vskip1mm

(ii)  If $\bU$ is  contained in $\sD (\bG)$, then it is unipotent. }
\vskip3mm

\ni{\it Proof.} Let $\bT$ be the unique maximal $k$-torus of the smooth connected nilpotent $k$-subgroup $\bU$. Then $\bT$ is normal, and hence central, in the connected group $\bB$. This implies that $\bB$ is contained in the centralizer  $\mathbf{M}:=Z_{\bG}(\mathbf{T})$ of $\bT$ in $\bG$. But then $\bB'\subseteq \mathbf{M}' = Z_{\bG'}(\mathbf{T}')$.   Note that $\mathbf{M}'$, being the centralizer of a torus in the connected reductive group $\bG'$, is a connected reductive subgroup. As some $\{\bB',\bB'\}$-double coset in $\bG'$ is open, some $\{\mathbf{M}',\mathbf{M}'\}$-double coset in the reductive group $\bG'$ is open. By Proposition 1.6 of [4] then $\mathbf{M}' = \bG'$, i.e., $\mathbf{T}'$ is central in $\bG'$.  Hence, the commutator subgroup $(\bG,\bT)_{\ok}$ is contained in the geometric unipotent radical $\sR_u(\bG_{\ok})$ of $\bG$. So the smooth connected $k$-subgroup $(\bG,\bT)$ is contained in the $k$-unipotent radical $\sR_{u,k}(\bG)$; i.e., $\bT$ commutes with $\bG$ modulo $\sR_{u,k}(\bG)$. Therefore (see Corollary 2 in [2], \S11.14), under the natural projection $\bG\rightarrow \bG/\sR_{u,k}(\bG)$, $Z_{\bG}(\bT)$ maps onto $\bG/\sR_{u,k}(\bG)$, hence $\bG = Z_{\bG}(\bT)\sR_{u,k}(\bG)$. But since $\sR_{u,k}(\bG)$ is a $k$-wound unipotent normal subgroup of $\bG$, it commutes with the torus $\bT$. This implies that $\bG = Z_{\bG}(\bT)$; i.e., $\bT$ is contained in the center of $\bG$.   
This proves the first assertion of the proposition.  
\vskip1mm

To prove the second assertion, assume that $\bU\subset \sD(\bG)$. Then $\bT'\subset \sD(\bG')$. But $\sD(\bG')$ is a connected semi-simple group, so it does not contain any nontrivial central tori. This implies that $\bT' = \pi(\bT_{\overline k})$ is trivial. As the kernel of $\pi$ is a unipotent group, we conclude that $\bT$ is trivial and hence $\bU$ is unipotent. This proves the second assertion.\hfill$\square$
\vskip2mm

\ni{\bf Proposition 2.} {\em Assume that either $\bG$ is perfect\,(and quasi-reductive) or it is pseudo-reductive. If the image $\bU'$ of $\bU_{\ok}$ in $\bG'$ is central, then $\bU$ is central in $\bG$ and $B$ is of finite index in $\sG$.}

\vskip3mm 

\ni{\it Proof.} Let us assume first that $\bG$ is perfect and quasi-reductive. Then $\bU$ is unipotent (Proposition 1({\it {ii}}). Now if the unipotent group $\bU'$ is central in the reductive group $\bG'$, then it is trivial. This implies that $\bU_{\ok}$ is contained in $\sR_u(\bG_{\ok})$ and hence $\bU\subseteq \sR_{u,k}(\bG)$. But in a perfect quasi-reductive group $\bG$, $\sR_{u,k}(\bG)$ is central (2.2), so $\bU$ is central. 

Now let us assume that $\bG$ is an arbitrary pseudo-reductive group and $\bU'$ is central in $\bG'$. Then the commutator subgroup $(\bG,\bU)_{\ok}$ is contained  in $\sR_u(\bG_{\ok})$, and so $(\bG,\bU)\subseteq \sR_{u,k}(\bG)$. But as $\bG$ has been assumed to be pseudo-reductive, its $k$-unipotent radical $\sR_{u,k}(\bG)$ is trivial, so $(\bG,\bU)$ is trivial and hence $\bU$ is central. 

If \,$\bU$ is central in $\bG$, then $\bU(k)$, which contains a subgroup of $U$ of finite index, is central in $G$, so $U$ is virtually central in $\sG$. Now from the 
Bruhat decomposition $\sG =\bigcup_{n\in N} BnB = \bigcup_{n\in N} UnB$, we see that if $U$ is virtually central (in $\sG$), then $B$ is of finite index in $\sG$. This proves the proposition. \hfill$\square$ 
\vskip2mm

\ni{\bf Proposition 3.} {\em  Assume that for every $s\in S$, the index of $B$ in $\sG_s$ is infinite. Let $\bV_{\ok}$ be a smooth connected unipotent $\ok$-subgroup of $\bG_{\ok}$ which is normalized by $B$. Let $V:=\bV_{\ok}(\ok)\cap \sG$ inside $\bG(\ok)$. Then $V$ is contained in $B$.}
\vskip3mm

\ni{\it Proof.}  We consider the subgroup $BV=VB$ of $\sG$. As it contains $B$, it equals  $\sG_X$ for a subset $X$ of $S$. To prove that $V$ is contained in $B$, we need to show that $X$ is empty.  Let us assume that this is not the case, and fix an $s\in X$. Let $V_s = V\cap \sG_s$. Then $BV_s =HUV_s= \sG_s$, and, in particular, $V_s$ is not contained in  $B$.  The pair $(B,N_s := N\cap \sG_s)$ is a Tits system in $\sG_s$ of rank 1. Let $\mathfrak B$ be the building associated with this Tits system. $\mathfrak B$ is just the infinite set  $\sG_s/B$ endowed with the left-multiplication action of $\sG_s$.  Any pair of distinct points in this building constitute an apartment. Let $x$ and $x'$ be the points of $\mathfrak B$ fixed by $B$ and $B' = sBs^{-1}$ respectively. Then $H\,(\subseteq B\cap sBs^{-1})$ fixes both $x$ and $x'$ and $s$ interchanges them.  As $\sG_s(= UV_sH)$ operates transitively on the set of apartments of $\mathfrak B$, we see that it acts $2$-transitively on $\mathfrak B$. Since $H$ fixes $x$ and $x'$, we conclude that  the subgroup $UV_s$ is 2-transitive on $\mathfrak B$. Now we recall that according to Proposition 1({\it i}), the unique maximal $k$-torus of the connected nilpotent subgroup $\bU$ is contained in the center of $\bG$.  Therefore, $\bU_{\ok}\bV_{\ok}$ is a nilpotent subgroup of $\bG_{\ok}$. Hence the subgroup $UV$, and so also the subgroup $UV_s$, is virtually nilpotent. By the following lemma, which was  pointed out to the author by Katrin Tent,  such a group cannot act 2-transitively on an infinite set. Thus we have arrived at a contradiction. This proves that $V$ is contained in $B$.  \hfill$\square$
\vskip2mm

\ni{\bf Lemma 1.} {\em A virtually nilpotent group cannot act $2$-transitively on an infinite set.}
\vskip2mm

\ni{\it Proof} (by Katrin\:Tent). We will prove the lemma by contradiction. Let $\sH$ be a virtually nilpotent group which acts $2$-transitively on an infinite set $\sX$. After replacing $\sH$ with its quotient by the kernel of the action, we may (and do) assume that the action is faithful. Then any nontrivial commutative normal subgroup of $\sH$ acts simply transitively on $\sX$. Let $\sN$ be a nilpotent normal subgroup of $\sH$ of finite index and $C(\sN)$ be its center.  Then $C(\sN)$ is a nontrivial commutative normal subgroup of $\sH$, so it acts simply transitively on $\sX$, in particular, $C(\sN)$ is infinite.  Fix an $x\in \sX$ and let $\sH_x$ be the stabilizer of $x$ in $\sH$. We identify $\sX$ with $C(\sN)$ using the bijection $g\mapsto gx$. With this identification, the action of $\sH_x$ on $\sX$ is just the conjugation action (of $\sH_x$) on $C(\sN)$. As $\sH_x$ acts transitively on $\sX-\{x\}$, we conclude that the conjugation action of $\sH_x$ on $C(\sN)$ is transitive on the set of nontrivial elements of  $C(\sN)$. But since $\sN$ is of finite index in $\sH$, the orbit of any element of $C(\sN)$ under $\sH$, and hence under $\sH_x$,  is finite.  This is a contradiction.   \hfill$\square$
\vskip2mm

For an arbitrary quasi-reductive group $\bG$, we have the following converse of Proposition 2.

\vskip2mm 
\ni{\bf Proposition 4.} {\em Assume that  the index of $B$ in $\sG$, as well as in $\sG_s$ for every $s\in S$, is infinite.  Then the image $\bU'$ of \,$\bU_{\ok}$ in $\bG'$  is not  central.} 
\vskip3mm

\ni{\it Proof.} Applying Proposition 3 to the smooth connected normal subgroup $\bV_{\ok} :=\sR_{u}(\bG_{\ok})$ of $\bG_{\ok}$ we see that $V:=\sR_{u}(\bG_{\ok})(\ok)\cap \sG$ is contained in  $B$. Note that  the kernel of the map $\pi:  \bG_{\ok}\rightarrow \bG' = \bG_{\ok}/\sR_u(\bG_{\ok})$ restricted to $\sG\subseteq G =\bG(k)\,(\subset \bG(\ok))$ is $V$. As the subgroup $\bnG_{\ok}$ maps onto  $\bG'$, and a subgroup of $\sG$ is Zariski-dense in $\bnG$, the image of $\sG$ in $\bG'(\ok)$ is Zariski-dense in $\bG'$. So  $UV/V$ is virtually central in $\sG/V$ (i.e., a subgroup of $UV/V$ of finite index is central in $\sG/V$) if and only if the image $\bU'$ of $\bU_{\ok}$ in $\bG'$ is central. Now assume, if possible, that $\bU'$ is central in $\bG'$, or, equivalently,  that $UV/V$ is virtually central in $\sG/V$, then, from the Bruhat decomposition $\sG =\bigcup_{w\in W^T}BwB = \bigcup_{w\in W^T}UwB= \bigcup_{w\in W^T}UVwB$, we see at once that $B/V$ is of finite index in $\sG/V$, and hence $B$ is of finite index in $\sG$, a contradiction.    \hfill$\square$ 
\vskip2mm

\ni{\bf Proposition 5.} {\em  Assume that the index of $B$ in $\sG$ is infinite and one of the following two conditions hold:  $(1)$ either $\bG$ is perfect\,(and quasi-reductive) or it is pseudo-reductive; 
$(2)$ $\bG$ is quasi-reductive and the index of $B$ in $\sG_s$ is infinite for every $s\in S$. Then the  commutator subgroup $(\bB,\bU)$ contains a nontrivial $k$-split smooth connected unipotent subgroup. In particular,  then $\sR_{us,k}(\bB)$ is nontrivial.}
\vskip2mm

\ni{\it Proof.} According to Proposition 2 and the preceding proposition, $\bU'$ cannot be central in $\bG'$.  Since $\bT'$ is central in $\bG'$ (Proposition 1({\it i})), but $\bU'$ is not, we conclude that $\bT'\ne \bU'$, and hence the unipotent radical $\sR_u(\bU')$ of $\bU'$ is nontrivial. If $\bU$ is  central in  $\bB$, then $\bU'$, and so also $\sR_u(\bU')$,  is central in $\bB'$. We assert that $\bB'$ cannot contain a nontrivial smooth unipotent central subgroup. For, otherwise, $\bB'(\ok)$ would contain a nontrivial unipotent element $z$ which is central (in $\bB'$).  The reduced centralizer $\bC'_z$ of $z$ in $\bG'$ would then contain $\bB'$.   As some $\{\bB',\bB'\}$-double coset  in $\bG'$ is open, some $\{\bC'_z,\bC'_z\}$-double coset is open too. But it has been shown by Martin Liebeck that for no nontrivial unipotent element $z\in \bG'(\ok)$, a $\{\bC'_z,\bC'_z\}$-double coset can be open in $\bG'$ ([8], Ch.\,1, Cor.\,8). This proves our assertion. So we conclude that $\bU$ cannot be central in $\bB$.
Then the commutator subgroup $(\bB, \bU)$ is a nontrivial smooth connected normal $k$-subgroup of $\bB$. We claim that this subgroup is unipotent. To verify this claim, we consider the maximal reductive quotient $\bB^{\rm{red}}:= \bB_{\ok}/\sR_u(\bB_{\ok})$  of $\bB_{\ok}$. The image of  the smooth connected nilpotent normal subgroup $(\bB,\bU)_{\ok}$ of $\bB_{\ok}$ in $\bB^{\rm{red}}$ is a smooth connected nilpotent normal subgroup contained in the (semi-simple) derived subgroup of $\bB^{\rm{red}}$, so it is trivial. This implies that $(\bB,\bU)_{\ok}\subseteq \sR_u(\bB_{\ok})$, and hence $(\bB,\bU)$ is unipotent, as claimed.   

We will prove that $(\bB,\bU)$ contains a nontrivial $k$-split smooth connected $k$-subgroup. Then the conjugates of this $k$-split subgroup of $(\bB,\bU)$ under $\bB(k_s)$ generate a smooth connected $k$-subgroup which, being contained in $(\bB,\bU)$, is unipotent; it is clearly  a normal subgroup of  $\bB$ and Theorem B.3.4 of [6] implies that it is $k$-split, hence $\sR_{us,k}(\bB)$ is nontrivial.  
\vskip1mm

As we will not need to work with the original $\bU$ anymore, for simplicity we will denote the smooth connected unipotent normal $k$-subgroup $(\bB,\bU)$ of $\bB$ by $\bU$ in the rest of this proof, and $\bU'$  will now  denote the image of $(\bB,\bU)_{\ok}$ in $\bG'$. Let us assume, if possible,  that $\bU$ does not contain any nontrivial $k$-split smooth connected subgroup. Then $\bU$ is $k$-wound, and hence every $k$-torus of  
$\bB$ commutes with it.  This implies that the subgroup $\bB_t$ of $\bB$ generated by $k$-tori commutes with $\bU$. But then every $\ok$-torus of $\bB_{\ok}$ commutes with  $\bU_{\ok}$ since according to Proposition A.2.11 of [6], $(\bB_t)_{\ok}$ is the subgroup of $\bB_{\ok}$ generated by all the $\ok$-tori in $\bB_{\ok}$.
\vskip1mm

\vskip1mm

We will now show that $\bB'$ contains a nontrivial smooth connected unipotent subgroup in its center, contradicting what we showed above, so this will prove Proposition 5.  Let us fix a Borel subgroup $\bS'\ltimes \bV'$ of $\bB'$, with $\bS'$ a maximal $\ok$-torus and $\bV'$ the unipotent radical of the Borel subgroup. Note that as $\bU'$ is a smooth connected unipotent normal subgroup of $\bB'$ it is contained in $\bV'$. Now we inductively define, for each positive integer $i$, the normal subgroups  $\bV'_i$ of  $\bV'$ contained in $\bU'$ as follows. Let $\bV'_1 = \bU'$, and having defined $\bV'_i$, let $\bV'_{i+1}$ be the commutator subgroup $(\bV', \bV'_i)$.  Let $n$ be the largest integer such that $\bV'_n$ is nontrivial. Then $\bV'_n$ is a nontrivial smooth connected subgroup of $\bU'$ which commutes with $\bV'$. Since every $\ok$-torus of $\bB_{\ok}$ commutes with $\bU_{\ok}$, the $\ok$-torus $\bS'$ of $\bB'$ commutes with $\bU'$, and hence also with the subgroup $\bV'_n$.     Thus the Borel subgroup 
$\bS'\ltimes \bV'$  of $\bB'$ commutes with the subgroup $\bV'_n$, which implies that $\bV'_n$ is central in $\bB'$ (by Corollary 11.4 of [2]).   \hfill $\square$
\vskip2mm

\ni{\bf Proposition 6.} {\em Under the hypothesis of the preceding proposition, there exists a proper pseudo-parabolic $k$-subgroup $\bP$ of $\bG$ such that (a) the $k$-split unipotent radical $\sR_{us,k}(\bP)$ of $\bP$ contains $\sR_{us,k}(\bB)$ and (b) $B\subseteq P:= \bP(k)$.}  
\vskip2mm

\ni{\it Proof.} According to the preceding proposition, $\sR_{us,k}(\bB)$ is nontrivial. Since $\sR_{u,k}(\bG)$ is $k$-wound, the image $\overline{\bV}$ of $\sR_{us,k}(\bB)$ in the pseudo-reductive quotient $\overline{\bG}:=\bG/\sR_{u,k}(\bG)$ is a nontrivial  $k$-split smooth connected  unipotent subgroup which is normalized by the image $\overline B$ of $B$ in $\overline {G} := \overline{\bG}(k)$. So by Theorem C.3.8 of [6], there exists a proper pseudo-parabolic $k$-subgroup $\overline{\bP}$ of $\overline{\bG}$ such that $\overline{\bV}\subseteq \sR_{us,k}({\overline{\bP}})$ and $\overline{B}\subseteq \overline{\bP}(k)$.  We can choose a $1$-parameter subgroup  $\lambda :{\rm{GL}}_1\rightarrow \bG$ such that for the induced  $1$-parameter subgroup $\overline{\lambda}: {\rm{GL}}_1\rightarrow \overline{\bG}$, $\overline{\bP} = P_{\overline{\bG}}(\overline{\lambda})$. Now take $\bP$ to be the pseudo-parabolic subgroup $P_{\bG}(\lambda)$ of $\bG$ (cf.\,2.4). This pseudo-parabolic subgroup clearly has the required properties. \hfill$\square$ 
\vskip2mm

\ni{\it Remark 4.} Proposition 1.6 of [4] and Corollary 8 of [8, Ch.\,1] used in the proof of Propositions 1 and 5 are implied by the following general result proved by Guralnick, Malle and Tiep independently and by a simple direct argument: {\it Let $\bG$ be a connected semi-simple group over an algebraically closed field $K$, $x$ a noncentral element of $\bG(K)$, $\bC_x$ the reduced centralizer of $x$ in $\bG$, then no $\{\bC_x,\bC_x\}$-double coset can be dense in $\bG$}.  See Corollary 5.5 in  [7].
\vskip2mm

\vskip5mm

{\centerline{\bf 4. Proof of Theorem A}}
\vskip3mm

We will continue to use the notation introduced in 3.1 and 3.2.

\vskip2mm

{\it Proof of Theorem A. } We begin by observing that since in Theorem A, it has been assumed that every $k$-split torus of $\bG$ is central (equivalently, $\sD(\bG)$ is $k$-anisotropic), $\bG$ does not contain a proper pseudo-parabolic $k$-subgroup, see 1.3.  But  according to Proposition 6, in case $B$ is of infinite index in $\sG$, and either $\bG$ is perfect\,(and quasi-reductive) or it is pseudo-reductive, $\bG$ does contain a proper pseudo-parabolic $k$-subgroup. Therefore, $B$ must be of finite index in $\sG$. 
\vskip1mm


\vskip1mm
We will assume now that either $\bG$ is perfect and quasi-reductive, or $\bG$ is reductive. In both the cases, $G$  (and hence $\sG$) is Zariski-dense in $\bG$. As $B$ is of finite index in $\sG$, $B$ is also Zariski-dense in $\bG$. Hence, the Zariski-closure $U\bU$ of $U$ (which is normalized by $B$)  is a nilpotent smooth normal $k$-subgroup of $\bG$. Its identity component $\bU$ is then a nilpotent smooth connected normal $k$-subgroup of $\bG$. If $\bG$ is reductive, such a subgroup is central.   On the other hand, if $\bG$ is perfect, by Proposition 1({\it{ii}}), $\bU$ is unipotent, so it is contained in the $k$-unipotent radical $\sR_{u,k}(\bG)$. But  if $\bG$ is perfect and quasi-reductive, $\sR_{u,k}(\bG)$ is central\,(2.2). Thus, in this case also,  $\bU$ is central. The centrality of $U\bU$ in $\bG$, and so of $U$ in $\sG$,  follows in both the cases now from Lemma 5.3.2 of [6] using the centrality of $\bU$ in $\bG$.

Now as $N$ normalizes $H$ and 
$U$ is central, $B =HU$ is normalized by $N$ and hence by $\sG$ which is generated by $B\cup N$. Therefore, $B$ is a normal subgroup of $\sG$. But the normalizer of  $B$ equals $B$ which implies that $B =\sG$. {\hfill $\square$}
\vskip2mm

\ni{\bf A Tits system of rank 1 in ${\mathrm{SL}}_1(D)$:} Let $\sf{F}$ be a field. Then $\sf{F}^{\times}$ acts on $\sf{F}$ by multiplication, and we form the semi-direct product $\mathscr{G} := \sf{F}\rtimes \sf{F}^{\times}$. The solvable group $\mathscr{G}$ admits a split Tits system $(B,N)$ of rank $1$ with $B =\{ (0,x)\vert\, x\in \sf{F}^{\times}\}$ and $N =\{(-1,-1), (0,1)\}$. 

Now let $D$ be the division algebra with center the $2$-adic field $\bQ_2$ and of dimension $d^2$ over $\bQ_2$. Let $\sf{F}$ be the field extension of degree $d$ of the field with $2$ elements. Let $\bG = {\rm{SL}}_{1,D}$. Then $\bG$ is an absolutely simple  simply connected algebraic group defined and anisotropic over $\bQ_2$, and $G := \bG(\bQ_2)$ is the subgroup ${\rm{SL}}_1(D)$ of $D^{\times}$ consisting of elements of reduced norm $1$. Let $G_2$ be the ``second congruence subgroup'' of $G$ as in [12], 1.1. Then $G/G_2 \simeq \mathscr{G} = \sf{F}\rtimes \sf{F}^{\times}$ ([12], \S1). Hence, $G$ admits a spherical Tits system of rank 1.  Such a Tits system cannot be weakly-split (cf.\,Theorem A).
\vskip5mm

\centerline{\bf 5. Proof of Theorem B when the Tits system is of rank 1}
\vskip3mm

\ni{\bf 5.1.}  Since Theorem B  is of special interest in case the Tits system is of rank 1, and the proof  is simpler than in the general case, we first prove it for Tits systems of rank 1. The proof of Theorem B for Tits systems of arbitrary rank will be given in the next section; it will use results proved in this section.  

We will assume in this section that the Tits system $(B,N)$ in $\sG$ is weakly-split, and of rank 1 (i.e., $S=\{s\}$) with, as before, $B =HU$, $H = B\cap N$ and $U$ is a nilpotent normal subgroup of $B$. We will also assume that  the index of $B$ in $\sG$ is infinite. To prove Theorem B (for Tits systems of rank 1), we need to prove that {(1) there exists a (proper) pseudo-parabolic $k$-subgroup $\bP$ of $\bG$ such that $B = \bP(k)\cap \sG$, and  (2) if the Tits system $(B,N)$ is saturated and $B$ does not contain a non-central normal subgroup of the group of $k$-rational points of a $k$-isotropic minimal perfect smooth connected normal $k$-subgroup  of $\mathbf{G}$, then the Tits system is a standard Tits system, i.e., it is as in 1.2. Hence the Tits system is split (possibly, in terms of a different  decomposition of $B$; see, however, Proposition 8 below).}
\vskip1mm

As the index of $B$ in $\sG$ has been assumed to be infinite, by Proposition 6, there is a proper pseudo-parabolic $k$-subgroup $\bP$ of $\bG$ such that  $B\subseteq P:= \bP(k)$.  Since $\bP$ is a proper pseudo-parabolic $k$-subgroup of  $\bG$, $P$ cannot contain $\sG$ (2.5).  Since the Tits system is of rank $1$, any subgroup of $\sG$ which properly contains $B$ equals $\sG$. Hence, $P\cap \sG= B$, and $\bP$ is a maximal proper pseudo-parabolic $k$-subgroup of $\bG$ (2.6). This proves assertion (1).  
\vskip2mm

The proof of assertion (2) requires Proposition 7 given below, whose proof in turn requires the following lemma. 
 We first fix some notation. 
Let $\Psi$ be an irreducible root system of rank $>1$,  given with a basis $\Delta$. Let $W$ be the Weyl group of $\Psi$ and $R =\{ r_a\ |\ a\in \Delta\}$, where for a root $a$, $r_a$ is the reflection in $a$; $R$ is a set of involutive generators of $W$.
We now  fix a root $a\in \Delta$. Let $\Psi'$ be the root subsystem of $\Psi$ spanned by $\Delta-\{a\}$ (i.e., $\Psi'$ is the intersection of the $\bQ$-span of $ \Delta-\{a\}$ with $\Psi$) and $W' (\subset W)$ be the Weyl group of $\Psi'$; $W'$ is generated by the subset $R-\{r_a\}$.  The length of an element $w\in W$, in terms of the generating set $R$ of $W$, will be denoted by $\ell(w)$. Let $w_0$ and $w_0'$ be the longest elements of $W$ and $W'$ respectively. Both these elements are of order $2$.
\vskip2mm
  
  If $\Psi$ is of type $A_m$ and $a$ is one of the end roots of $\Delta$ then $\# W'\backslash W/W' = 2$,  as can be easily verified (see the example in 5.3 below). On the other hand, we have the following:
  \vskip2mm

\ni{\bf Lemma 2.}  {\em Let us assume that the root system $\Psi$ is of rank $m>1$. Then $\#W'\backslash W/W'>2$ unless the root system is of type $A_m$ and $a$ is one of the two end roots of $\Delta$.}
\vskip2mm

\ni{\it Proof.} If the root system is not reduced (then it is of type $BC_m$), the subset $\Psi^{\bullet}$ consisting of nondivisible roots of $\Psi$ is a root system of type $B_m$ with Weyl group $W$, and $\Delta$ is a basis of $\Psi^{\bullet}$.  Therefore, to prove the lemma, we may (and do) replace $\Psi$ by $\Psi^{\bullet}$ and assume that $\Psi$ is reduced (and irreducible). We will denote by $\Psi^+$ the positive system of roots determined by the basis $\Delta$.
\vskip1mm

Let us first consider the case where $w_0 = -1$ (this is the case unless the root system $\Psi$  is of type $A_m\, (m>1)$, $D_m$, with $m$ odd, or $E_6$) and assume for the sake of contradiction that $\#W'\backslash W/W' \leqslant 2$.  Then,  clearly,  $W = W'\cup W'r_aW'$, and we conclude that $w_0 = -1= w_1r_aw_2$, with $w_1,\,w_2\in W'$.  Then $-r_a = w_1^{-1}{w_2}^{-1}$. As $a$ is the only positive root  which is transformed into a negative root by $r_a$, we see that $w_1^{-1}{w_2}^{-1}\,(\in W')$ takes all the roots in ${\Psi'}^+:=\Psi'\cap \Psi^+$ into negative roots. Therefore, $(w_2w_1)^{-1} = w_0'$. Hence, $w_2= w_0'w_1^{-1}$, and so $\ell(w_2) = \ell(w_0')-\ell(w_1)$. But then $\ell(w_0) = \ell(w_1r_aw_2) \leqslant \ell(w_1)+1+\ell(w_2)= \ell(w_0') +1$. Since $2\ell(w_0)=\#\Psi$ and $2\ell(w_0') =\#\Psi'$, we obtain the bound $\#\Psi\leqslant \#\Psi' +2$. But it is easily seen that this bound does not hold. 
\vskip1mm

The following argument to prove the lemma in case $w_0\ne -1$, and more generally if the root system $\Psi$ is simply laced,  was kindly provided by John Stembridge.  
\vskip1mm   

We will assume now that the root system $\Psi$ is simply laced, is of rank $m>1$, and if it is of type $A_m$ then $a$ is not an end root of $\Delta$. If $a$ is not an end root (i.e., if it is connected to at least 2 simple roots) we set $i = 1$, $a_1 = a$, and pick any  two simple roots  to which $a$ is connected and call them $a_2$ and $a_3$. On the other hand, if $a$ is an end root (then $\Psi$ is one of the following types: $D_m$, $m\geqslant 4$, $E_6$, $E_7$, $E_8$)  let $b$ be the unique simple root which is connected to three other simple roots. Let $a_1= a$, $a_2$, $\ldots$\,, $a_i=b$ be the nodes on the shortest path joining $a$ to $b$ in the Dynkin graph of $\Delta$, and let $a_{i+1}$ and $a_{i+2}$ be the two simple roots, other than $a_{i-1}$, connected to $b = a_i$.  For $j\leqslant i+2$, let $r_j= r_{a_j}$. Then   $r_1r_2$, $r_2r_3$, $\ldots$\,, $r_{i-1}r_i$, 
$r_ir_{i+1}$, and $r_ir_{i+2}$ all have order 3. 
\vskip1mm

Let $w = r_1r_2\cdots r_i\cdot r_{i+1}r_{i+2}\cdot r_i\cdots r_2r_1 \in W$. By considering the monoid of words in $S$ modulo substitutions according to the braid relations and $s^2 = 1$, for $s\in S$, it follows from Proposition 5 in \S1.5 of Ch.\,IV of [3] that any expression in $S$ for an element of $W$ is brought to a reduced expression via such substitutions and that any two reduced expressions are converted into each other exclusively by substitution using the braid relations.  Thus, the given expression for $w$ is reduced and the only other reduced expression for $w$ is the one obtained by interchanging the commuting $r_{i+1}$ and $r_{i+2}$ (note that $a_{i+1}$ cannot be connected to $a_{i+2}$ since both are connected to $a_i$).   Hence, $w$ is of length $2i+2$ and it has exactly two reduced expressions, both of them begin and end with $r_1 = r_a$.  We now claim that the double coset $W' w W'$ is distinct from $W'$ and $W' r_1 W'$ (so $\#  W'\backslash W/W' >2$).  If $w$ belongs to $W'$ then there is reduced expression for $w$ which is free of $r_1$, and if it belongs to $W'r_1W'$, then it has a reduced expression in which either the first or the last term is different from $r_1$, a contradiction.\hfill$\square$ 
\vskip2mm

In the next proposition, which will be used again in the next section, we assume that the Tits system $(B,N)$ in $\sG$  is weakly-split and of rank 1, but we do {\it not} assume that it is saturated. Let  $\bP$ be a pseudo-parabolic $k$-subgroup of $\bG$, $P = \bP(k)$, and assume that $B=P\cap \sG$. Then we have the following:
\vskip2mm

\ni{\bf Proposition 7.} {\em The pseudo-parabolic $k$-subgroup $\bP$ contains all $k$-isotropic minimal perfect smooth connected normal $k$-subgroups of $\bG$ except one. The one not contained in $\bP$ is of $k$-rank $1$.}
\vskip2mm

\ni{\it Proof.} As in 2.3 let $\bnG$ be the identity component of the Zariski-closure of $G$ in $\bG$. It is a quasi-reductive  $k$-subgroup of $\bG$, $\bnP = \bP\cap \bnG$ is a pseudo-parabolic $k$-subgroup of $\bnG$, and $\snG := \bnG(k)\cap \sG$ is of finite index in $\sG$.  Any perfect smooth connected $k$-subgroup of $\bG$ is contained in $\bnG$ (2.3) and every $k$-torus and every $k$-split smooth connected unipotent $k$-subgroup of $\bP$ is contained in $\bnP$ (2.4). So, in particular, $\bS\subset \bnP$.  
\vskip1mm

We  fix a minimal pseudo-parabolic $k$-subgroup $\bP_0$ contained in $\bP$ and a maximal $k$-split torus $\bS$ of $\bG$ contained in $\bP_0$. Then $\bP_0 = Z_{\bG}(\bS)\ltimes \sR_{us,k}(\bP_0)$ (Proposition C.2.4 and Corollary 2.2.5 of [6]). We will denote the $k$-Weyl group $N_{\bG}(\bS)(k)/Z_{\bG}(\bS)(k)$ of the quasi-reductive $k$-group $\bG$ by $W$. Let $P_0 = \bP_0(k)$,  $Z=Z_{\bG}(\bS)(k)$, $\sP_0=P_0\cap \sG$ and $\sN = N_{\bG}(\bS)(k)\cap \sG$. Since $\sR_{us,k}(\bP_0)(k)\subset G^+\subseteq \sG$, $P_0 = \sP_0Z$. Hence, for all $n\in \sN$, $P_0nP_0 =  \sP_0 n\sP_0Z$.  We recall from 1.2 that $\sN$ maps onto $W$. So we see using the Bruhat decomposition of $G$ with respect to $P_0$ (that is $G=\bigcup_{n\in \sN}P_0nP_0$, Theorem C.2.8 of [6]) that $G = \sG Z$. Hence, $P = (P\cap \sG)Z = BZ$. 
\vskip1mm

Let $\Phi$ be the set of $k$-roots of $\bG$ with respect to the maximal $k$-split torus $\bS$, $\Phi^+ (\subset \Phi)$ be the positive system of roots and $\Delta$ the set of simple roots determined by the  minimal pseudo-parabolic $k$-subgroup $\bP_0$.  Let $R = \{ r_a\ | \ a\in  \Delta\}$, where for $a\in \Delta$, $r_a\in W$ is the reflection in $a$.   Then $R$ is a set of involutive generators of $W$; the lengths of elements in $W$ will be in terms of the generating set $R$. Let $w_0$ be the longest element of $W$. For a subset $Z$ of $R$, we will denote by $W_Z$ the subgroup of $W$ generated by the elements in $Z$.  Let $X$ be the subset of $R$ such that $P = P_0 W_X P_0$. Then since $\bP$ is a maximal proper pseudo-parabolic $k$-subgroup of $\bG$, $X = R-\{r_a\}$, for an $a\in \Delta$. Let $\Delta_a$ be the connected component of $\Delta$ containing $a$, and $\Psi$ the irreducible component of $\Phi$ spanned by $\Delta_a$. Since $\bP\supseteq \bP_0\supset Z_{\bG}(\bS)$ and the set of nonzero weights of $\bS$ on ${\rm{Lie}}(\bP)$ 
contains all the connected components of $\Phi$ except $\Psi$, it is clear that $\bP$ contains all $k$-isotropic minimal perfect smooth connected  normal $k$-subgroup of $\bG$ except the one whose root system with respect to $\bS$ is $\Psi$ (see Proposition C.2.32 of [6]). Let $Y$ be the set of reflections in the roots contained in $\Delta_a$, and $Y' =Y-\{r_a\}$; $W_Y$ is clearly the Weyl group of the root system $\Psi$. Since $\Delta_a$ is orthogonal to $\Delta-\Delta_a$, $W = W_Y\times W_{R-Y}$. Now as  $\sG=B\cup BsB$ and $B\subseteq P$, $\sG$ (so also $G=\sG Z$, as $Z\subset P_0\subseteq P$)  is contained in the union of  two $\{ P, P\}$-double cosets. From this we infer using the Bruhat decomposition (Theorem C.2.8 of [6]) that $W$ is the union of two $\{W_X,W_X\}$-double cosets, and hence $W_Y$ is  the union of two $\{W_{Y'},W_{Y'}\}$-cosets.  (Therefore, $W_Y = W_{Y'}\cup W_{Y'}r_aW_{Y'}$.)

\vskip1mm

We shall now show that $\Psi$ is of rank $1$ using the fact that the Tits system under consideration is weakly-split.  Assume, if possible, that $\Psi$ is of rank $m>1$. As $W_Y$ is  the union of two $\{W_{Y'},W_{Y'}\}$-cosets, from Lemma 2 we conclude that $\Psi$  is of type $A_m$ (so $\Psi$ is reduced),  and $a$ is one of the two end roots of the basis $\Delta_a$ of this root system. Since $W =W_Y\times W_{R-Y}$ and $W_Y = W_{Y'}\cup W_{Y'}r_aW_{Y'}$, we obtain that $W = W_X\cup W_Xr_aW_X$, and then $w_0$ will have to belong to $W_Xr_aW_X$ which implies that $W_Xr_aW_X = W_Xw_0W_X$.  Hence,  $G = P\cup Pw_0P$. As the Weyl group of  the Tits system $(B=P\cap \sG,N)$ is $\{1,s\}$, $H\subseteq B\cap sBs^{-1} \subseteq P\cap gPg^{-1}$ for some $g\in Pw_0P$. Therefore, $H$ is contained in $p(P\cap w_0Pw_0^{-1})p^{-1}$ for some $p\in P$. Let $p = xz$, with $x\in B$ and $z\in Z$. Then 
$p(P\cap w_0Pw_0^{-1})p^{-1}=x(P\cap w_0Pw_0^{-1})x^{-1}$.  We will  show below that if $\Psi$ is the root system of type $A_m$ with $m>1$ and $a$ is an end root of its basis $\Delta_a$, then there cannot exist a nilpotent normal subgroup $U$ of $B=P\cap \sG$ and a subgroup $H$ of $B$ which is contained in $x(P\cap w_0Pw_0^{-1})x^{-1}$, for some $x\in B$, such that $B = HU$. This will prove the proposition. 

\vskip1mm
Suppose there exist  such $H$ and $U$. Let $\bU$ be the identity component of the Zariski-closure of $U$. Then as $\bnP$ is the identity component of the Zariski-closure of $P\cap\snG$ in $\bG$, and so also that of $B=P\cap\sG$, cf.\,2.5,  and $U$ is a nilpotent normal subgroup of $B$, we see that $\bU$ is a nilpotent normal subgroup of $\bnP$. Now we note that by Proposition  3.5.12\,(1) of [6], $\bnQ:= \bnP\cap w_0\bnP w_0^{-1}$ is a smooth connected $k$-subgroup. It clearly contains $\bS$. Since $H$ is contained in $x(P\cap w_0Pw_0^{-1})x^{-1}$ for some  $x\in B$ and $\bnP$ is the identity component of the Zariski-closure of $P$ in $\bG$ (2.5), the identity component of the Zariski-closure of $H$ in $\bG$ is contained in $x\bnQ x^{-1}$ for an $x\in B$. From this we obtain, using again the fact that  $\bnP$ is the identity component of the Zariski-closure of $B=HU$ in $\bG$, that $\bnP = x\bnQ x^{-1}\bU$. Since $B$ normalizes both $\bnP$ and $\bU$ and $x$ lies in $B$, we conclude that $\bnP = \bnQ\bU$. As the root system $\Psi$ is reduced, we now see from Proposition 3.3.5 of [6] that  if an element of $\Psi$ is a weight of $\bS$ on 
${\rm{Lie}}(\bnP)$, then it is a weight on either  ${\rm{Lie}}(\bnQ)$ or ${\rm{Lie}}(\bU)$.

\vskip1mm

To determine the weights of $\bS$ on ${\rm{Lie}}(\bnP)$ and ${\rm{Lie}}(\bnQ)$ contained in $\Psi$, we enumerate the roots in $\Delta_a$ as $\{a_1, a_2, \ldots, a_m\}$ so that $a_1 = a$ and for $i\leqslant m-1$, $a_i$ is {\it not} orthogonal to $a_{i+1}$, i.e., the corresponding nodes are connected in the Coxeter graph.   Then $w_0(a_i) = -a_{m+1-i}$ for all $i\leqslant m$. The  weights of $\bS$ on ${\rm{Lie}}(\bnP)$ which belong to $\Psi$ are all the roots in 
$\Psi^+:= \Psi\cap \Phi^+$, and also all the negative roots contained in the span of $\{a_2,\dots, a_m\}$. The weights of $\bS$ on ${\rm{Lie}}(w_0\bnP w_0^{-1})$ contained in $\Psi$ are all the roots in $-\Psi^+$, and also all the positive roots contained in the span of $\{a_1,\ldots, a_{m-1}\}$. From this we see that the weights of $\bS$ on ${\rm{Lie}}(\bnP)$ which lie in $\Psi$, but are not weights on  ${\rm{Lie}}(\bnQ)$, are the $m$ roots $a_i+\cdots +a_m$, with $i\geqslant 1$.   These $m$ roots must therefore be weights of $\bS$ on ${\rm{Lie}}(\bU)$. Note that for $i\geqslant 2$, $-(a_i+\cdots +a_m)$ is a weight on  ${\rm{Lie}}(\bnQ)$. Now for a root $b\in \Psi$, let $\bU_b$ be the corresponding root group (see C.2.20 of [6]); $\bU_b$ is the largest smooth connected unipotent $k$-subgroup of $\bnG$  normalized by $\bS$ on whose Lie algebra the only weight of $\bS$ is $b$.  If $b \in \Psi$ is a weight of $\bS$ on  ${\rm{Lie}}(\bnQ)$, then 
$\bU_b$ is contained in both $\bnP$ and $w_0\bnP w_0^{-1}$ and so it is contained in $\bnQ=\bnP\cap w_0\bnP w_0^{-1}$. On the other hand, if $b \in \Psi$ is a nonzero weight of $\bS$ on ${\rm{Lie}}(\bnP)$ which is not a weight on  ${\rm{Lie}}(\bnQ)$, then it must be a weight on  ${\rm{Lie}}(\bU)$ and the root group $\bU_b$  is contained in $\bU$.  Let now $b = a_2+\cdots+a_m$. Then $-b$ is a weight of $\bS$ on  ${\rm{Lie}}(\bnQ)$, so $\bU_{-b}\subset \bnQ$, and $b$ is not a weight of $\bS$ on  ${\rm{Lie}}(\bnQ)$ so it is a weight on ${\rm{Lie}}(\bU)$ and $\bU_{b}\subset \bU$. Since $\bU$ is a (connected nilpotent) normal subgroup of $\bnP$, it is normalized by $\bU_{-b}$, and hence $\bU_{-b}\bU$ is a solvable subgroup ($\bU_{-b}\bU$ is in fact  nilpotent  since according to Proposition\:1({\it{i}}) the unique maximal torus of the nilpotent group $\bU$ is contained in the center of $\bG$). But the subgroup $\bU_{-b}\bU$ contains the subgroup $\bnG_b$ generated by 
$ \bU_{\pm b}$ which is a nonsolvable group (nonsolvability of $\bnG_b$ can be seen by considering the image of $(\bnG_b)_{\ok}$ in the reductive group $\bG'= \bG_{\ok}/\sR_u(\bG_{\ok})$; this image is isomorphic to either ${\rm{SL}}_2$ or ${\rm{PGL}}_2$). This is impossible, so we are done. 
\hfill $\square$
\vskip2mm

\ni{\bf 5.2.} To prove assertion (2) of 5.1, we will assume now that  the Tits system $(B,N)$ is saturated and $B$ does not contain a non-central normal subgroup of the group of $k$-rational points of a $k$-isotropic minimal perfect smooth connected normal $k$-subgroup  of $\mathbf{G}$. Then it follows from the first observation in Remark 1, Proposition 7, the description of pseudo-parabolic $k$-subgroups of $\bG$  and Proposition C.2.32 of [6]  that $\bP$ does not contain any $k$-isotropic perfect smooth connected normal $k$-subgroups of $\bG$ and the  $k$-rank of $\sD(\bG)$ is 1. Hence $\bP$ is a minimal pseudo-parabolic $k$-subgroup of $\bG$. It is known that for $g\in G$, $\bP\cap g{\bP}g^{-1}$ contains the centralizer of a maximal $k$-split torus of $\bG$ (Proposition C.2.7 of [6]); moreover, since $\sD(\bG)$ is of $k$-rank $1$, we easily see that any two distinct minimal pseudo-parabolic $k$-subgroups of $\bG$ are opposed to each other. Therefore, for any $g\notin P$, as $gPg^{-1}\ne P$ (recall that $P$ equals its normalizer in $G$), $g\bP g^{-1}$ is different from $\bP$, and so it is opposed to $\bP$. Hence,  $\bP\cap g{\bP}g^{-1}$ equals the centralizer (in $\bG$) of  some maximal $k$-split torus of $\bG$. Then as the Tits system $(B,N)$ is saturated, $H = P\cap sPs^{-1}\cap \sG = Z_{\bG}(\bS)(k)\cap\sG$ for a suitable maximal $k$-split torus $\bS$ of $\bG$. 
From the fact that the identity component of the Zariski-closure of $H= Z_{\bG}(\bS)(k)\cap \sG$ in $\bG$ equals $Z_{\bnG}(\bS)$ (2.5), we see that the Zariski-closure of $H$ in $\bG$ contains $\bS$ as the maximal $k$-split central torus, and hence the normalizer of $H$ in $G$ is contained in $N_{\bG}(\bS)(k)$. Now since $N$ normalizes $H$ and contains it as a subgroup of index $2$, whereas  $\sN :=N_{\bG}(\bS)(k)\cap \sG$ is the unique such subgroup of $\sG$, we conclude that $N =  \sN$. {\it Thus the saturated Tits system $(B,N)$  is a  standard Tits system.}  This proves assertion (2) of 5.1.

\vskip2mm
\ni {\bf 5.3. A Tits system of rank 1 in $G:={\rm{SL}}_{n+1}(k)$}. Let $B$ be the subgroup of matrices in  $G$ whose first column has all the entries  zero except the top entry. Let $B'$ be the subgroup of matrices in $G$ whose last column has all the entries zero except the bottom entry. Let $H = B\cap B'$, and choose a $g\in G$ such that $gBg^{-1} = B'$, and let $N = H\cup gH$. Then $N$ is a subgroup and $H$ is a normal subgroup of $N$ of index $2$. It is easily checked that $(B,N)$ is a spherical Tits system in $G$ of rank $1$.
Note that $B$ itself admits a Tits system of rank $n-1$ since it is isomorphic to ${\rm{GL}}_{n}(k)\ltimes k^n$.

It has been pointed out by Pierre-Emmanuel Caprace and Katrin Tent that giving a Tits system $(B,N)$ of rank 1 in $G$ is equivalent to providing a $2$-transitive action of $G$ on a set $X$ with at least two elements: Given $X$ with a $2$-transitive action of $G$, and a pair $x\ne x'$ in $X$, let $B$ be the stabilizer of $x$ and $N$ be the stabilizer of the subset $\{x,x'\}$. Then $(B,N)$ is a Tits system of rank $1$ in $G$. Conversely, given a Tits system $(B,N)$ of rank $1$ in $G$, the natural action of $G$ on $G/B$ is $2$-transitive. Now in the above example, the set $X$ is the projective space ${\bP}^n(k)$ with the natural action of ${\rm{SL}}_{n+1}(k)$. 
\vskip5mm

\centerline{\bf 6. Proof of Theorem B for Tits systems of arbitrary rank}
\vskip3mm

\ni {\bf 6.1.} We will continue to use the  notation introduced in \S\S 1-2. Let $(B,N)$ be a weakly-split Tits system in $\sG$ with Weyl group $W^T$. We assume that $W^T$ is finite, i.e., the Tits system is spherical, and let $S$ be  the distinguished set of involutive generators of $W^T$.  Let $B = HU$, with $H = B\cap N$, and $U$ a nilpotent normal subgroup of $B$. Let $\bU$ be the identity component of the Zariski-closure of $U$ in $\bG$. For $s\in S$, let $\sG_s = B \cup BsB$. For a subset $X$ of $S$, let $W^T_X$ be the subgroup of $W^T$ generated by the elements in $X$, and $\sG_X = BW^T_XB$. Then for each $s\in S$, and $X\subseteq S$, $\sG_s$ and $\sG_X$ are subgroups of $\sG$.  

\vskip2mm

\ni{\bf 6.2.} We will prove Theorem B by induction on the rank ($=\#S$) of the Tits system $(B,N)$. If the rank is  zero, i.e., if $S$ is empty, $B =\sG$. Then the first assertion of Theorem B holds if we take $\bP =\bG$. The hypothesis  ``$B$ does not contain a non-central normal subgroup of the group of $k$-rational points of a $k$-isotropic minimal perfect smooth connected normal $k$-subgroup  of $\mathbf{G}$''  in the second assertion  does not hold in this case if $\sD(\bG)$ is $k$-isotropic.  So to prove the second assertion when $S$ is empty, we assume that $\sD(\bG)$ is $k$-anisotropic. If the Tits system $(B,N)$, with $B=\sG$ is saturated, then $H :=B\cap N= \bigcap_{n\in N}n\sG n^{-1} = \sG$, and hence, $N =\sG$. The Tits system $(B,N) = (\sG,\sG)$ is clearly the unique standard Tits system in $\sG$ if $\sD(\bG)$ is $k$-anisotropic.  
\vskip1mm

 Now let us assume that $\#S>0$. For $s\in S$, as the index of $B$ in $\sG_s\,(\subseteq \sG)$ has been assumed to be  infinite, the index of $B$ in $\sG$ is infinite. Then according to Proposition 5, $\sR_{us,k}(\bB)$ is nontrivial, and by Proposition 6 there exists  a proper pseudo-parabolic $k$-subgroup $\bQ$ of $\bG$ such that $\sR_{us,k}(\bB)\subseteq \sR_{us,k}(\bQ)$ and $B\subseteq \bQ(k)$. If  $B =\bQ(k)\cap \sG$, then the first assertion holds if we take $\bP = \bQ$, so let us assume that $B\ne \bQ(k)\cap \sG$. As $\bQ(k)\cap\sG$ properly contains $B$, and is not equal to $\sG$\,(2.5),  it equals $\sG_X$ for a nonempty subset $X$\,$(\ne S)$ of $S$. The quotient $\overline{\bM} := \bQ/\sR_{us,k}(\bQ)$ is a quasi-reductive $k$-group.  Let $\osM$ be the image of $\sG_X\subseteq \bQ(k)$ under the natural projection $\bQ(k)\rightarrow \overline{\bM}(k)$ and  let $K= \sR_{us,k}(\bQ)(k)$. As $\sR_{us,k}(\bQ)$ is a $k$-split smooth connected unipotent subgroup, $K\subset G^+\subseteq \sG$. So $K$ is a normal subgroup of $\sG_X=\bQ(k)\cap\sG$, and the projection map $\sG_X\rightarrow \osM$ induces an isomorphism of $\sG_X/K$ with $\osM$. We will use this isomorphism to identify $\sG_X/K$ with $\osM$.  We see from Proposition 3 that $K$ is contained in $B$. 
\vskip1mm

\ni{\bf 6.3.} Let $N_X = N\cap \sG_X$, $\overline{B} = B/K$ and $\overline{N}_X = N_XK/K$.  Then $(B,N_X)$ is a Tits system in $\sG_X$ and $(\overline{B}, \overline{N}_X)$ is a Tits system in $\osM=\sG_X/K$. Let 
$\overline{H}$ (resp., $\overline{U}$) be the image of $H$ (resp., $U$) in $\overline{B}$. As $B =HU$, $\overline{B} = \overline{H} \overline{U}$; moreover, since $B\cap N_XK= (B\cap N_X)K =HK$, 
$\overline{B}\cap \overline{N}_X=\overline{H}$, so the Tits system $(\overline{B}, \overline{N}_X)$ in $\osM$ is weakly-split. As $H\subseteq N_X\cap HK\subseteq N_X\cap B\subseteq H$, the natural homomorphism $W^T_X = N_X/H\rightarrow \overline{N}_X/\overline{H}$ is an isomorphism, and hence the Weyl group of the Tits system $(\overline{B},\overline{N}_X)$ in $\osM$ is $W^T_X$.  For $x\in X$, $\osM_x \, (=\overline{B}\cup \overline{B}x\overline{B})$ is the image of $\sG_x$ in $\osM$ and the induced map $\sG_x/B\rightarrow \osM_x/\overline{B}$ is bijective. Therefore, the index of $\overline{B}$ in $\osM_x$ is infinite for all $x\in X$. Now since the rank of the Tits system $(\overline{B}, \overline{N}_X)$ is $\#X< \#S$,  we conclude by induction on the rank of Tits systems that there exists a proper pseudo-parabolic $k$-subgroup  $\overline{\bP}$ of the quasi-reductive $k$-group $\overline{\bM} = \bQ/{{\sR}_{us,k}(\bQ)}$ such that $\overline{B} = {\overline{\bP}}(k)\cap \osM$. Let $\bP$ be the inverse image of $\overline{\bP}$ in $\bQ$. Then $\bP$ is a pseudo-parabolic $k$-subgroup of $\bG$ (Lemma 3.5.5 of [6]), and $B\subseteq \bP(k)$. Moreover, since $B$ contains the kernel $K = \sR_{us,k}(\bQ)(k)$ of the natural map $\bQ(k)\rightarrow \overline{\bM}(k)$,  and its image 
$\overline{B}$ in $\overline{\bM}(k)$ equals $\overline{\bP}(k)\cap\osM$, we conclude that $B = \bP(k)\cap\sG =:\sP$.  This proves the first assertion of Theorem B.
\vskip1mm

\ni{\bf 6.4.} Now to prove the second assertion  of Theorem B, we will assume that the Tits system $(B,N)$ is saturated and $B$ does not contain a non-central normal subgroup of the group of $k$-rational points of a $k$-isotropic minimal perfect smooth connected normal $k$-subgroup  of $\mathbf{G}$. We wish to first find a maximal $k$-split torus $\bS$ of $\bG$ such that $\sZ:= Z_{\bG}(\bS)(k)\cap \sG\subseteq H = B\cap N$. For this purpose, and also for use later, we introduce the following notation: for an element $w\in W^T$, $\ell(w)$ will denote its length with respect to the distinguished set of generators $S$ of the Weyl group $W^T$. The longest element of $W^T$ will be denoted by $w_0$.  We will now show that $B\cap w_0Bw_0^{-1} = H$.  For this, we will use Lemma 13.13 of [13]. This lemma says that $B\cap ww'B(ww')^{-1}\subseteq B\cap wBw^{-1}$ if $\ell(ww') = \ell(w)+\ell(w')$. But given any $w \in W^T$, let $w'= w^{-1}w_0$. Then $ww' = w_0$, and $\ell(w) +\ell(w') = \ell(w_0)$, which implies that $H =\bigcap_{w\in W^T} wBw^{-1} =B\cap w_0Bw_0^{-1}$.  Therefore, $H = B\cap w_0Bw_0^{-1}= \sP\cap w_0\sP w_0^{-1}$. Now we note that $\bP\cap w_0\bP w_0^{-1}$ is a smooth connected subgroup which contains the centralizer $Z_{\bG}(\bS)$ of a maximal $k$-split torus $\bS$ of $\bG$ (Propositions 3.5.12 and  C.2.7 of [6]). Hence, $\sZ  \subseteq \sP\cap w_0\sP w_0^{-1}= H$. 
\vskip1mm

\ni{\bf 6.5.} We will next prove that $\bP$ is a minimal pseudo-parabolic $k$-subgroup of $\bG$. Let us fix a minimal pseudo-parabolic $k$-subgroup $\bP_0$ of $\bG$ contained in $\bP$ and containing $Z_{\bG}(\bS)$. Let $P_0 = \bP_0(k)$, $\sP_0 = P_0\cap \sG$, $Z = Z_{\bG}(\bS)(k)$ and $\sN = N_{\bG}(\bS)(k)\cap \sG$. Recall from 1.2 that $\sN$ maps onto the $k$-Weyl group $W = N_{\bG}(\bS)(k)/Z_{\bG}(\bS)(k)$. It was observed in the proof of Proposition 7 that for all $n\in \sN$, $P_0nP_0 = \sP_0n\sP_0Z$ and $G = \sG Z$.

Let $\Phi$ be the set of $k$-roots of $\bG$ with respect to $\bS$, $\Phi^+$ the positive system of roots and $\Delta\,(\subseteq \Phi^+)$ the set of simple roots given by $\bP_0$. For a nondivisible root $b\in \Phi$, let $\bU_b$ be the root  group corresponding to it (see [6], C.2.20); $\bU_b$ is the largest smooth connected unipotent $k$-subgroup of $\bG$ normalized by $\bS$ on whose Lie algebra the only weights of $\bS$ under the adjoint action are positive integral multiples of $b$.  Then there is a subset $\Delta'$ of  $\Delta$ such that $\bP$ is generated by $\bP_0$ and the root groups $\bU_{-a}$, $a\in \Delta'$. We will now prove that $\bP=\bP_0$,  or, equivalently, that $\Delta'$ is empty. The condition that $B\,(=\sP=\bP(k)\cap \sG)$ does not contain a non-central subgroup of the group of $k$-rational points of a $k$-isotropic  minimal perfect smooth connected normal $k$-subgroup of $\bG$ is equivalent to the condition that $\Delta'$ does not contain any connected components of $\Delta$ (see Remark 1 above and Proposition C.2.32 of [6]).  \vskip1mm

Let  $R = \{r_a\,|\, a\in \Delta\}$, where for $a\in \Delta$, $r_a\in W$ is the reflection in $a$. The reflections in $R$ generate $W$; for $w\in W$, we will denote the length of $w$, with respect to this generating set, by $\ell(w)$. The rank of $W$, i.e., the number of elements in $\Delta$,  is equal to the $k$-rank  of $\sD(\bG)$\,(1.2).   For a subset $X$ of $R$, we will denote by $W_X$ the subgroup of $W$ generated by the elements in $X$. For $w\in W$, $\sP_0 w \sP_0$ and $\sP w \sP$ will denote the double cosets $\sP_0 n \sP_0$ and $\sP n \sP$ respectively, where $n$ is any representative of $w$ in $\sN$. For $X\subseteq R$, let $\sG^X = \sP_0 W_X\sP_0$. Any subgroup of $\sG$ that contains $\sP_0$ equals $\sG^X$ for a unique $X\subseteq R$.  Now let $Y\subset R$ be such that $B=\sP = \sG^Y$.  Then $Y =\{ r_a\ |\ a\in \Delta'\}$.  For $s\in S$, as $\sG_s = B\cup BsB$ is a subgroup containing $B =\sP$, there is a  subset $Y_s$ of $R$ containing $Y$ such that $\sG_s = \sG^{Y_s}$. Since there is no subgroup of $\sG$  lying properly between $B$ and $\sG_s$,  $Y_s = Y\cup \{r_{a_s}\}$ for some $a_s\in \Delta-\Delta'$, and then $BsB = \sP r_{a_s}\sP$.   For distinct elements $s,s'$ of $S$, since $\sG_s\cap \sG_{s'} = B =\sP$,  and the collection of subgroups $\{\sG_s\}_{s\in S}$ generates $\sG$, we conclude that if $s\ne s'$ then $a_s\ne a_{s'}$, and moreover, $\{ a_s\,|\,s\in S\} = \Delta-\Delta'$.
\vskip1mm

Assume, if possible, that $\Delta'$ is nonempty. Then since $\Delta'$ does not contain any connected components of $\Delta$, we can find a root $a\in \Delta'$ such that there is a root $b\in \Delta-\Delta'$ connected to $a$. As $b\notin \Delta'$, $-b$ is not a weight of $\bS$ on the Lie algebra of $\bP$.  Let $s\in S$ be the element such that $b=a_s$. Let $\bP'$ be the pseudo-parabolic $k$-subgroup generated by $\bP$ and the root group $\bU_{-b}$ and let $P' =\bP'(k)$. Let ${\overline{\bM}}' =\bP'/\sR_{us,k}(\bP')$ be the maximal quasi-reductive quotient of $\bP'$ and ${\overline{M}}' :={\overline{\bM}}'(k)$. Then $\sP':=P'\cap\sG = \sG_s = \sP\cup \sP r_b\sP$.  Since $\bP'\supset \bP$, $\sR_{us,k}(\bP')\subset \sR_{us,k}(\bP)$ (Proposition 3.5.14 of [6]) and hence $K':=\sR_{us,k}(\bP')(k)\subset \sP=B$. Let $\overline{\sP}$, $\overline{\sP}'$  and $\overline{N}'$ denote the respective images of $\sP$, $\sP'$ and $N' := N \cap \sP'$ in $\overline{M}'=P'/K'$.  Since $(B=\sP,N')$ is a Tits system of rank 1 in $\sP'$, $(\overline{\sP}, \overline{N}')$ is a Tits system of rank 1 in $\overline{\sP}'$. 
To see that this Tits system is weakly-split, let $\oH$ and $\oU$ be the images of $H$ and $U$ in $\overline{\sP}$. Then $\overline{\sP}= \oH\oU$; moreover, as $\sP\cap N'K' = (\sP\cap N')K' \subseteq (\sP\cap N)K' = (B\cap N)K' =HK'$, we find that ${\overline{\sP}}\cap {\overline{N}}' = \oH$, which shows that the Tits system $({\overline{\sP}},{\overline{N}}')$ is weakly-split.

By Proposition C.2.32 of [6], there is a $k$-isotropic minimal perfect smooth connected normal $k$-subgroup of $\overline{\bM}'$ whose $k$-root system is the connected component of the root system (of $\overline{\bM}'$)  containing $a$. As $b$ is connected to $a$, this connected component contains $b$ and so is of rank at least 2.  Proposition 7 applied to the quasi-reductive group $\overline{\bM}'$, the image of $\bP$ in $\overline{\bM}'$ (this image is a pseudo-parabolic $k$-subgroup of $\overline{\bM}'$) in place of $\bG$ and $\bP$, and the weakly-split Tits system $(\overline{\sP}, \overline{N}')$ of rank  $1$ in $\overline{\sP}'\,(\subseteq \overline{M}'=\overline{\bM}'(k))$, implies  that the image of $\bP$ in $\overline{\bM}'$ contains  this $k$-isotropic minimal perfect smooth connected normal $k$-subgroup of $\overline{\bM}'$. From this we infer that $-b$ is a weight of $\bS$ on the Lie algebra of $\bP$. But this is not true. Thus we have proved that $\Delta'$ is empty, and hence $\bP =\bP_0$ is a minimal pseudo-parabolic $k$-subgroup of $\bG$. In particular, the map $s\mapsto r_{a_s}$ from $S$ into $R$ is bijective.  

\vskip1mm

\ni{\bf 6.6.} Having proved that $\bP$ is a minimal pseudo-parabolic $k$-subgroup of $\bG$ and $B =\sP =\bP(k)\cap \sG$, we will now prove that $H= \sZ =Z\cap \sG$ and $N =N_{\bG}(\bS)(k)\cap \sG=N_{\sG}(Z)$.  We begin by recalling from [3], Ch.\,IV, \S 2.4, Theorem 2,  that 
\vskip2mm

$(*)$\ \  {\it for $w\in W^T$ (resp., $w\in W$) and $s\in S$ (resp., $r\in R$) $\ell(sw)>\ell(w)$ (resp., $\ell(rw) >\ell(w)$) if and only if $BsB\cdot BwB = BswB$ (resp., $\sP r \sP \cdot \sP w\sP = \sP rw \sP$), and if \,$\ell(sw)<\ell(w)$ then $BsB\cdot BwB$ is the union of two distinct $\{ B,B\}$-double cosets $BwB$ and $BswB$.}
\vskip1mm

\vskip1mm

For $s\in S$, let $a_s$ be as in 6.5. We will first prove that the map $S\rightarrow W$ which maps $s$ onto $r_{a_s}$ extends to an isomorphism of $W^T$ onto $W$.   For this purpose, let $s,s'\in S$ and $m$ be the common order of $ss'$ and $s's$ in $W^T$. Let $w$ (resp., $\widetilde{w})$ be the product of the first $m$ terms of the sequence $\{s', s, s', s,\ldots \}$ (resp., of the sequence $\{r_{a_{s'}}, r_{a_s}, r_{a_{s'}}, r_{a_s},\ldots \}$). Then $\ell(w)=m$, and $\ell(sw)=m-1$. Since $B =\sP$, $BsB = \sP r_{a_s}\sP$,  and $Bs'B = \sP r_{a_{s'}}\sP$, by repeated application of $(*)$ we easily see that $\ell({\widetilde{w}}) = m$, and $\ell(r_{a_s}{\widetilde{w}}) = m-1$. This observation implies that the order of  the product $r_{a_{s'}}r_{a_s}$ (in the $k$-Weyl group $W$) is also $m$. Now since $(W^T,S)$ and $(W,R)$ are Coxeter groups, we conclude that the map $s\mapsto r_{a_s}$ of $S$ into $R$ extends uniquely to an isomorphism $W^T\rightarrow W$. In particular,  $W^T$ and $W$ are of equal order, say $n$. 
\vskip1mm

We have shown in 6.4 that $H\supseteq \sZ$. We will now prove that these two subgroups are equal. For this, we note that since $\bP$ is a minimal pseudo-parabolic $k$-subgroup of $\bG$, $B= \sP$ has precisely $n$ distinct 
conjugates in $\sG$ containing $\sZ$, these are $w\sP w^{-1}, w\in W$. The intersection of these conjugates equals $\sZ$. On the other hand, since the Weyl group $W^T$ of the Tits system $(B,N)$ is of order $n$, at least $n$ distinct conjugates, namely $wBw^{-1}$, $w\in W^T$,  of $B$ contain $H$, and hence also $\sZ$.  Therefore, there are exactly $n$ distinct conjugates of $B=\sP$ which contain $H$, and these conjugates are $wBw^{-1}, w\in W^T$.  Thus the collections  $\{w\sP w^{-1} \  |\ w\in W\}$ and   $\{wBw^{-1}\ |\ w\in W^T\}$ are the same. The hypothesis that the Tits system $(B,N)$ is saturated now gives that  $H = \bigcap_{w\in W^T}wBw^{-1} = \bigcap_{w\in W}w\sP w^{-1} = \sZ$. 
\vskip1mm

To complete the proof of the second assertion of Theorem B it only remains to show that $N=N_{\bG}(\bS)(k)\cap \sG = \sN$. For this, we recall from 2.5 that the identity component of the Zariski-closure of $H=\sZ =Z_{\bG}(\bS)(k)\cap \sG$ equals $Z_{\bnG}(\bS)$, and hence $\bS$ is the maximal $k$-split central torus of the Zariski-closure of $H$ in $\bG$. This implies that any element of $G$ which normalizes $H$ belongs to $N_{\bG}(\bS)(k)$. Therefore, $N\subseteq N_{\sG}(H)\subseteq N_{\bG}(\bS)(k)\cap \sG = \sN$.  But $[\sN:H] = n = [N:H]$ as both the Weyl groups $W$ and $W^T$ have order $n$. This implies that $N =\sN$, and we have proved  Theorem B.\hfill$\square$

\vskip2mm

We will end this paper with the following proposition and two remarks.
\vskip2mm

{\ni}{\bf Proposition 8.} {\em Assume that $\bG$ is a perfect pseudo-reductive $k$-group, $\sG$ is a subgroup of \,$G:=\bG(k)$ that is Zariski-dense in $\bG$ and contains $G^+$. Let $(B,N)$ be a split spherical Tits system in $\sG$ and $U$ be a nilpotent normal subgroup of $B$ such that $B=H\ltimes U$, with $H = B\cap N$. Assume that $B$ does not contain a non-central normal subgroup of the group of $k$-rational points of a $k$-isotropic minimal perfect smooth connected normal $k$-subgroup  of $\mathbf{G}$. Let $\bP$ be the minimal pseudo-parabolic $k$-subgroup of $\bG$ such that $B =\bP(k)\cap \sG$. Then $U = \sR_{u,k}(\bP)(k)$.}  
\vskip2mm

\ni{\it Proof.}  Let $\bH$ and $\bU$ be respectively the identity components of the Zariski-closures of $H$ and $U$ in $\bG$. Then $B$ normalizes $\bU$ and the identity component of the Zariski-closure of $B$ in $\bG$  is $\bH\bU$. As $\bG$ has been assumed to be perfect, $G$ is Zariski-dense in $\bG$ so $\bnG = \bG$, $\bnP =\bP$,  and $B =\bP(k)\cap\sG$ is Zariski-dense in $\bP$\,(2.3--2.5). Hence, $\bH\bU = \bP$. By Proposition 1({\it ii}), $\bU$ is unipotent. Since it is normalized by $B$, it is a smooth connected unipotent normal $k$-subgroup of $\bP$. This implies that $\bU\subseteq \sR_{u,k}(\bP)$. Note that $\sR_{u,k}(\bP) =\sR_{us,k}(\bP)$ by Corollary 2.2.5 of [6] since $\bG$ has been assumed to be pseudo-reductive. Let $\bS$ be a maximal $k$-split torus of $\bP$ such that $H = Z_{\bG}(\bS)(k)\cap \sG$\,(6.6). Then $\bH \subseteq \bM:= Z_{\bG}(\bS)$.  Let  $M = \bM(k)$, $P =\bP(k)$, and $S =\bS(k)$.   As $\bP$ is a minimal pseudo-parabolic $k$-subgroup of $\bG$, $\bP = \bM \ltimes \sR_{us,k}(\bP)$ (Proposition C.2.4 and Corollary 2.2.5 of [6]), hence $P = M\ltimes \sR_{us,k}(\bP)(k)$.   

Now  $\bP = \bH\bU \subseteq \bM \bU\subseteq \bM\ltimes \sR_{us,k}(\bP)=\bP$, which implies that $\bH = \bM$ and $\bU =\sR_{us,k}(\bP)$.   So $\bU(k)\subset G^+\subseteq \sG$. To prove that $U = \bU(k)$, we consider the natural projection map $\pi : P=M\ltimes \bU(k)\rightarrow M$. Its kernel is $\bU(k)$ which contains a subgroup of finite index of $U$, so $F:=\pi(U)$ is a finite subgroup of $M$ and $U\subseteq F \cdot \bU(k)$.   On the other hand, since $B = \bP(k)\cap \sG\supset \bU(k)$ and $B = H\ltimes U$,  we see that $\bU(k) \subset F\cdot U$. From this we easily deduce that the index of  $V:=U\cap \bU(k)$ in $\bU(k)$ is $\# F$, in particular this index is finite. The subgroup $V$ is clearly normalized by $S\cap B = S\cap \sG$.

We assert  that  $S\cap G^+$, and hence also $S\cap \sG$, is Zariski-dense in $\bS$. To prove this, we may clearly replace $\bG$ with a $k$-split semi-simple subgroup which contains $\bS$ as a maximal $k$-torus and whose root system  is the set of non-multipliable roots in the root system of $\bG$ with respect to $\bS$ (for a proof of existence of such a $k$-split semi-simple subgroup, see Theorem C.2.30 of [6]), and assume that $\bG$ is a $k$-split semi-simple group. But for a $k$-split semi-simple group  our assertion is well-known. (It can be proved by considering the simply connected central cover of $\bG$.)

We will show that  $\bU(k)$ does not contain a proper subgroup of finite index that is normalized by a subgroup $\sS$  of $S$ which is Zariski-dense in $\bS$.  Assuming this for a moment,  the subgroup $V := U\cap \bU(k)$,   which is a subgroup of $\bU(k)$ of finite index ($=\#F$) normalized by the group $S\cap \sG$ (which is Zariski-dense in $\bS$), equals $\bU(k)$. Hence, $\# F=1$, which implies that  $U\subseteq \bU(k)$. This in turn implies that $U=\bU(k)$, which proves the proposition modulo the assertion that $\bU(k)$ does not contain a  proper subgroup of finite index which is normalized by a subgroup $\sS$ of $S$ which is Zariski-dense in $\bS$. To prove that $\bU(k)$ does not contain a proper subgroup of finite index that is normalized by such an $\sS$, we note that as $\sS$ is Zariski-dense in $\bS$, for any root $a$, $a(\sS)$ is Zariski-dense in ${\mathrm{GL}}_1$. So the subfield $k_a$ of $k$ generated by $a(\sS)$  is infinite. Now, if $a$ is not multipliable, then the root group $\bU_a$ is a vector group admitting a $k$-linear structure such that any  element $s$ of $S$ acts by  multiplication by $ a(s)$. Thus any subgroup of $\bU_a(k)$ which is stable under the conjugation action of $\sS$ is a $k_a$-vector subspace of $\bU_a(k)$. Therefore, any $\sS$-stable subgroup of $\bU_a(k)$ either equals $\bU_a(k)$ or is of infinite index in it.  On the other hand, if $a$ is multipliable, then both $\bU_a/\bU_{2a}$ and $\bU_{2a}$ are nonzero vector groups and a similar consideration implies that any $\sS$-stable subgroup of $\bU_a(k)$ either equals $\bU_a(k)$ or is of infinite index in it.  Since $\bU(k)$ is generated by the $\bU_a(k)$'s, we are done.      
\hfill{$\square$}

\vskip3mm

\ni{\it Remark 5.} The above proposition does not  hold in general  if $\bG$ is not perfect. To see an example, consider the group $G={\mathrm{GL}}_2(\bR)$ with its standard Tits system $(B,N)$ in which  $B$ is the subgroup of upper triangular matrices and $N$ is the normalizer of the diagonal subgroup. The subgroup $H:=B\cap N$ is the subgroup of diagonal matrices. Let $U$ be the commutative normal subgroup of $B$ consisting of  the matrices
$ \begin{pmatrix} e^t & te^t  \\ 0 & e^t \end{pmatrix}, \  t\in \bR$.  Then $B = H\ltimes U$ is a nonstandard splitting of $B$.  The Zariski-closure of $U$ in ${\mathrm{GL}}_2$ is a connected subgroup which is not unipotent (it contains the central ${\mathrm{GL}}_1$) which shows that Proposition 1({{\it ii}) does not hold  in general if we drop the hypothesis that $\bU\subset \sD(\bG)$.
\vskip2mm

\ni{\it Remark 6.} Theorem B has the following simple consequences. Let $k_1$ and $k_2$ be infinite fields. For $i =1 ,\, 2$, let $\bG_i$ be a perfect quasi-reductive $k_i$-group. Assume that $\bG_i$ does not contain a $k_i$-isotropic perfect smooth connected proper normal $k_i$-subgroup.  Let $G_i^+$ be the normal subgroup of $\bG_i(k_i)$ generated by the group of $k_i$-rational points of the $k_i$-split unipotent radicals  of pseudo-parabolic $k_i$-subgroups of $\bG_i$ and $\sG_i$ be a subgroup of $\bG_i(k_i)$ containing $G_i^+$ that is Zariski-dense in $\bG_i$. If there exists an isomorphism $f:\,\sG_1\rightarrow \sG_2$ of abstract groups, then $k_1$-rank $\bG_1$ = $k_2$-rank $\bG_2$, and given a maximal $k_1$-split torus $\bS_1$ in $\bG_1$, there is a maximal $k_2$-split torus $\bS_2$ of $\bG_2$ such that $f$ carries the centralizer $Z_{\sG_1}(\bS_1)$ of $\bS_1$ in $\sG_1$ isomorphically onto the centralizer $Z_{\sG_2}(\bS_2)$ of $\bS_2$ in $\sG_2$.  In fact,  in $\sG_i$ we have the standard split Tits system provided by the underlying algebraic $k_i$-group $\bG_i$\,(1.2). Theorem B implies that under the isomorphism $f: \sG_1\rightarrow \sG_2$, these Tits systems are conjugate to each other. This at once implies both the assertions.  

Now for $i=1,\, 2$, let $\bG_i$ be a perfect pseudo-reductive $k_i$-group, and assume that $\bG_i$ does not contain a $k_i$-isotropic perfect smooth connected proper normal $k_i$-subgroup. We will use  the notation introduced in the preceding paragraph. From Proposition 8 and Theorem B we see that the  isomorphism $f:\,\sG_1\rightarrow \sG_2$ carries the  group of $k_1$-rational points of the $k_1$-split unipotent radical of a minimal pseudo-parabolic $k_1$-subgroup of $\bG_1$ isomorphically onto the group of $k_2$-rational points of the $k_2$-split unipotent radical of a minimal pseudo-parabolic $k_2$-subgroup of $\bG_2$.  
\vskip1mm

We note that Propositions 7.3.3(2) and 9.7.4(1) of [6] show that  over any local or global function field $k$ of characteristic 2 or 3, there exists a non-reductive pseudo-split perfect pseudo-reductive $k$-group $\mathbf{G}_1$ and a split connected semi-simple $k$-group $\mathbf{G}_2$ such that $\mathbf{G}_1(k) \simeq \mathbf{G}_2(k)$.

\vskip2mm

\vskip3mm

\noindent{\bf Acknowledgements.} I am grateful to Martin Liebeck for proving a  result used in the proof of  Proposition 5 above. His proof will appear in [8]. I thank Pierre-Emmanuel Caprace and Timoth\'ee Marquis for their paper [5] which inspired this paper.  I thank  Katrin Tent  for the proof of Lemma 1 and for several useful conversations and correspondence, and thank John Stembridge for the proof of Lemma 2. I thank Brian Conrad for carefully reading several versions of this paper and  for his comments  which have led to considerable improvement in the exposition. I thank the referee  for her/his comments and corrections.
\vskip1mm

I was partially supported by the NSF (grant DMS-1001748). The hospitality and support of the Institute for Advanced Study and the IHES are gratefully acknowledged. 

\vskip5mm

{\centerline{\bf References}}

\vskip3mm

\noindent 1. P.\,Abramenko,  M.\,Zaremsky, {\it Some reductive anisotropic groups that admit no non-trivial split spherical $BN$-pairs}, a paper posted on arXiv.
\vskip1mm

\noindent 2.  A.\,Borel, {\it Linear Algebraic Groups}, 2nd edition, Springer-Verlag, New York (1991).
\vskip1mm
 
\noindent 3. N.\,Bourbaki, {\it Lie algebras and Lie Groups}, Chapters IV, V, VI, Springer-Verlag, New York (2002).
 \vskip1mm
 
\noindent 4. J.\,Brundan, {\it Double coset density in exceptional algebraic groups}, J.\,London Math.\,Soc.\,{\bf 58} (1998), 63-83.
\vskip1mm

\noindent 5. P-E.\,Caprace, T.\,Marquis, {\it Can an anisotropic reductive group admit a Tits system?}, Pure and Appl.\,Math.\,Q.\,{\bf 7}(2011), 539-557. 
\vskip1mm

\noindent 6. B.\,Conrad, O.\,Gabber, G.\,Prasad, {\it Pseudo-reductive groups}, 2nd edition, New Mathematical Monographs \#17, Cambridge U.\:Press (to appear in 2013).

\vskip1mm

\noindent 7. R.\,Guralnick, G.\,Malle, P.H.\,Tiep, {\it Products of conjugacy classes in finite and algebraic simple groups}, preprint (February 2012).
\vskip1mm

\noindent 8. M.\,Liebeck, G.\,M.\,Seitz,  {\it Unipotent and nilpotent classes in simple algebraic groups and Lie algebras}, AMS Monograph (2012). 
\vskip1mm

\noindent 9. K.\,Tent, H.\,Van Maldeghem, {\it Moufang polygons and irreducible spherical $BN$-pairs of rank $2$}, I, Adv.\,Math. {\bf 174} (2003), 254-265.

\vskip1mm

\noindent 10. K.\,Tent, {\it Moufang polygons and irreducible spherical $BN$-pairs of rank $2$, the octagons}, Adv.\,Math. {\bf 181}(2004), 308-320.
\vskip1mm

\noindent 11. T.\,De Medts, F.\,Haot, K.\,Tent, H.\,Van Maldeghem, {\it Split $BN$-pairs of rank at least $2$ and the uniqueness of splittings}, J.\,Group Theory {\bf 8}(2005), 1-10.
\vskip1mm

\noindent 12. G.\,Prasad, M.\,S.\,Raghunathan, {\it Topological central extensions of SL$_1$(D)}, Inventiones Math.\,{\bf 92}
(1988), 645-689. 
\vskip1mm

\noindent 13. J.\,Tits, {\it Buildings of spherical type and finite $BN$-pairs}, Lecture Notes in Mathematics No.\:386, Springer-Verlag, New York (1974). 
\vskip1mm

\noindent 14. J.\,Tits, R.\,M.\,Weiss, {\it Moufang Polygons}, Springer Monographs in Mathematics, Springer-Verlag, New York (2002).

\vskip5mm





\noindent{Department of Mathematics}

\noindent{University of Michigan}

\noindent{Ann Arbor   MI 48109}

\noindent{e-mail: gprasad@umich.edu}
\end{document}